\documentclass[12pt]{article}

\usepackage{easytable,graphics,graphicx,amsmath,amssymb,array,eqnarray,verbatim,color,setspace,gensymb,authblk,comment,rotating,abstract,ulem,subfigure,amsthm,multirow,booktabs,stmaryrd}
\usepackage[version=4]{mhchem}
\usepackage{times}
\usepackage{multicol,enumitem,adjustbox,blindtext,cite,hyperref,xcolor}

\usepackage{soul}

\def\jacv#1{{\color{blue} #1}}

\def\mt#1{{\color{green} #1}}

\def\bypass{\textnormal{per,bypass}}
\def\bypassret{\textnormal{out,bypass}}
\def\feed{\textnormal{in}}
\def\exit{\textnormal{out}}
\def\per{\textnormal{per}}
\def\perrecy{\textnormal{recy}}
\def\ret{\textnormal{out}}
\def\lobnd{\textnormal{lo}}
\def\upbnd{\textnormal{up}}
\def\resp{\textit{resp. }}
\def\ie{\textit{i.e.},}

\def\permeance{\mathit{PM}}

\newcounter{eqnind}
\newcounter{eqnalind}

\newenvironment{eqnal}
   {\setcounter{eqnalind}{\theequation}
     \setcounter{equation}{\theeqnind}%
     \renewcommand\theequation{W\arabic{equation}}
     \align}
   {\endalign\setcounter{eqnind}{\value{equation}}\setcounter{equation}{\theeqnalind}}
\allowdisplaybreaks

\newcolumntype{C}[1]{>{\centering\arraybackslash}p{#1}}

\usepackage[margin=0.7in]{geometry}

\begin{document}

% ========== Edit your name here

\title{\textbf{Optimal Design of Membrane Cascades for Gaseous and Liquid Mixtures via MINLP}}
\date{}

\author[1]{Jose Adrian Chavez Velasco $^\ddagger$}
\author[1]{Radhakrishna Tumbalam Gooty $^\ddagger$}
\author[2]{\authorcr Mohit Tawarmalani$^\dagger$}
\author[1]{Rakesh Agrawal$^*$}
\affil[1]{\normalsize Davidson School of Chemical Engineering, Purdue University, West Lafayette, IN-47906}
\affil[2]{Krannert School of Management, Purdue University, West Lafayette, IN-47906 \vspace{-6ex}}

\date{}
\maketitle
\noindent Corresponding authors:  \texttt{mtawarma@purdue.edu}$^\dagger$,
\texttt{agrawalr@purdue.edu}$^*$\\
\noindent $^\ddagger$ These authors contributed equally 
\vspace{6pt}
\hrule
\begin{abstract}
\normalsize\noindent

Given the growing concern of reducing CO\textsubscript{2} emissions, it is desirable to identify, for a given separation carried out through a membrane cascade, the optimum design that yields the lowest energy consumption. Nevertheless, designing a membrane cascade is challenging since, there are often multiple feasible configurations that differ in their energy consumption and cost. In this work, we develop a Mixed Integer Non-linear Program (MINLP) that, for a given binary separation, which may be either liquid or gaseous, finds the cascade and its operating conditions that minimize energy consumption. To model the separation at each membrane in the cascade, we utilize the analytical solution of a system of differential and algebraic equations derived from the crossflow model and the solution-diffusion theory. We provide numerical evidence which shows that our single-stage membrane model accurately predicts experimental data. Unfortunately, the resulting membrane model is non-convex and, even state-of-the-art solvers struggle to prove global optimality of the cascades and the operating conditions identified. In this paper, we derive various cuts that help with relaxation quality and, consequently, accelerate convergence of branch-and-bound based solvers. More specifically, we demonstrate, on various examples, that our cuts help branch-and-bound solvers converge within 5\% optimality gap in a reasonable amount of time and such a tolerance level was not achieved by a simple formulation of the membrane model. The proposed optimization model is an easy-to-use tool for practitioners and researchers to design energy efficient membrane cascades.\\

\end{abstract}
\hrule

\pagebreak

\section{Introduction}
Membrane technology has gained considerable interest in recent years. Distinctive features such as modular scale-up flexibility, operational simplicity and relatively low capital cost make them attractive for the separation of a variety of mixtures. A few applications that have seen commercial success include water desalination at large scale, lube oil dewaxing, reduction of \ce{CO2} in natural gas, distributed scale nitrogen production from air, and so forth \cite{Galizia2017,WHITE200626}.  Besides these examples, a spectrum of applications exists where membranes can potentially be used on a large scale (see \cite{Galizia2017,WHITE200626}). The development of advanced materials with enhanced properties, such as higher area-to-volume ratio, has contributed substantially towards expanding the potential use of membranes. Despite the progress, membranes still have at most moderate permeabilities and selectivity  values in a majority of applications. Consequently, multiple membrane modules connected in a sequence, referred to as a \textit{membrane cascade} or simply a \textit{cascade} hereafter, is needed for separations which require one of the components to be recovered at high purity and high recovery. Nevertheless, identifying the most attractive cascade, even for binary separations, remains a challenging task.\\

It is well-known that, for a given separation, multiple membrane cascades exist. These cascades differ in terms of the overall energy input and cost \cite{SPILLMAN1995589,AGRAWAL19971029,QI1998209}. This energy requirement and cost also depend on the operating conditions of the cascade. The resulting energy consumption of separation is often a large fraction of the overall energy input of the chemical plant and this energy is usually sourced from the combustion of fossil fuels. Consequently, the separation may contribute significantly to the carbon footprint of the plant. Given rising concerns about greenhouse emissions, it is essential that energy efficient cascades be designed. Not surprisingly, this problem has been studied extensively in the literature. Nevertheless, we show that the current literature has not addressed global optimization of high-fidelity models. This is particularly important because most models for the permeation process use nonlinear nonconvex equations, and the local optimization techniques and/or meta-heuristics used in the literature do not guarantee that the discovered cascades are globally optimal, or even close in energy consumption to such a cascade. Here, we provide evidence that sometimes the local optimization solvers, or simply \textit{local solvers} hereafter, converge to a suboptimal solution that requires much higher energy than the global optimal solution. Besides, our computational experience suggests that, without a good starting solution, the local solvers fail to converge in a majority of cases. In addition to finding the most efficient membrane cascade, the global optimization solution is also important for comparing the membrane technology against an alternate separation technology. Although global optimization techniques are available for some separation technologies, such as distillation \cite{Caballero2006,Zheyu2018,Nallasivam2016,TUMBALAMGOOTY201913}, they are not yet available for design of membrane cascades making it difficult to perform head-to-head comparisons. Besides, since membrane material research is an active research endeavor, global optimization techniques are needed to ascertain selectivity and permeabilities of components at which membrane technologies may out-compete other separation processes, thereby helping to identify separations for which new membrane cascades are likely to yield most promising results. These reasons motivate the development of an effective global optimization approach for designing membrane cascades. \\

Towards the goal, we present several modeling advances in this article.  The highlights of the article are summarized in the following. First, we describe a unified \textit{permeator model} that is applicable for both gaseous and liquid mixtures. In its default form, the permeator model is a differential-algebraic equation (DAE) system. We express it as a system of algebraic equations by solving the DAE system analytically. This is in contrast to the common practice in the literature, where the DAE system is solved approximately using various discretization techniques. The analytical solution addresses the trade-off between the accuracy and the complexity of the model, which we shall discuss in detail in section 3 (\textsection 3). Second, we present a novel unified mixed-integer nonlinear program (MINLP) that is formulated to identify the membrane cascade which minimizes the overall power input. We solve the MINLP using an off-the-shelf global optimization solver, or simply a \textit{global solver} hereafter, such as BARON \cite{tawarmalani2005polyhedral,Kilin2018}. A global solver has built-in heuristics that generate multiple good-quality initial points for local search. Further, it continues local search with new initial points as it explores various parts of the search region until a proof of optimality is obtained. Once a global solver converges, the incumbent solution is guaranteed to be within a specified tolerance from the global optimum.  Third, for all of our test cases, global solvers fail to solve the MINLP to the desired optimality tolerance in a reasonable amount of time. To address this challenge, we derive additional cuts for the problem using physical insights and by exploiting the mathematical properties of the governing equations. Through numerical examples, we demonstrate that the additional cuts expedite the convergence characteristics of BARON and we are able to solve the MINLP within 5\%-optimality tolerance. Fourth, we apply our methods to two industrially important applications: separation of propylene/propane mixture and the separation of p-xylene from a mixture of xylene isomers.\\

The rest of the paper is organized as follows. We present a survey of the existing literature in \jacv{\textsection 2}. We describe the unified permeator model for a single stage in \jacv{\textsection 3}, and present our MINLP formulation along with computational experiments in \jacv{\textsection 4}. We discuss two case studies in \jacv{\textsection 5}, and present concluding remarks in \jacv{\textsection 6}.

\section{Literature review}
There are two aspects to designing a cascade for a given separation. First, identify the \textit{candidate set}; which we define as the set of all potentially attractive cascades for a given separation. Second, from the candidate set, identify the cascade along with its operating condition that optimizes the desired objective. The approaches used to generate the candidate set can be broadly classified as follows: (i) heuristic/intuition/empirical observation based approach, and (ii) superstructure based approach. In the first approach, the candidate set is identified based on physical insights and engineering judgements (see \cite{Qiu1989,BHIDE199113,BHIDE1993209,XU1996115,AGRAWAL1996129,LABABIDI1996185,QI1998209,HAO2008108,AHMAD2012119}). In the second approach, a superstructure is postulated, and the candidate set is obtained by discarding appropriate units and/or connections (nodes and/or arcs) from the superstructure (see \cite{AgrawalXu1996,AGRAWAL19971029,PATHARE2010263,Aliaga2017,QI199871,QI20002719,Uppaluri2004,UPPALURI2006832,SCHOLZ20151,RAMIREZSANTOS2018346,ARIAS2016371,ADI2016379,KUNDE2018164,MARRIOTT20034991}). We will adopt the second approach in this work.\\

% While their model is expressed in a general form and allows the use of any flux model for the transport process within the membrane, the authors, for their case studies, used flux equations derived using the solution-diffusion theory. 

Next, a permeator model is needed to assess a membrane cascade. Marriott and S\o rensen\cite{MARRIOTT20034975} used detailed mass, momentum, and energy balances to model the transport process on both retentate and permeate sides. Their model also takes into account concentration polarization and non-ideal solution properties. The resulting model is a system of partial differential algebraic equations. While this model accurately describes the process, it is computationally challenging to use it within an optimization framework. As such, most models in the literature make various assumptions that simplify the equations describing the transport process on both sides of the membrane. For example, the \textit{crossflow model} of Weller and Steiner \cite{weller1950separation,weller1950erratum}, when combined with flux equations derived from the solution-diffusion theory, yields such a model that is used in the literature. Here, the derivation assumes that the bulk phase concentration, velocity, and temperature gradients in the transverse direction (the direction perpendicular to the surface of the membrane) are negligible on both sides of the membrane. Consequently, the permeator model reduces to a simple ordinary differential algebraic equation (DAE) system. We use this simplified permeator model in our current work. Despite the simplifying assumptions, it has been shown that the model agrees well with experimental membrane performance data \cite{Pan1983,AHSAN201647}. In this work, we provide additional evidence that validates the cross-flow model against experimental data. One reason for a good agreement of the cross-flow model is that most commercial membranes are either asymmetric membranes or a membrane with a thin dense layer on a porous support. The works, \cite{Feng1999,QI199611,Qi1997} further accounted for pressure drop along the flow direction. However, as we will show in \textsection3, even without this change, the model agrees well with experimental data. Therefore, we will neglect the effect of pressure drop along the flow direction. Aliaga et al. \cite{Aliaga2017} and Adi et al. \cite{ADI2016379} additionally assume that the retentate and permeate sides are perfectly mixed. This leads to the well-known \textit{perfect mixing model} \cite{weller1950separation}, which can be described using a system of algebraic equations. However, as we will show in \textsection 3, the predictions from the perfect mixing model are less reliable than those from the crossflow model. Besides assuming perfect mixing, Aliaga et al. \cite{Aliaga2017} also assume that the separation factor (defined as $(y^\per/(1-y^\per))/(x^\feed/(1-x^\feed)$, where $y^\per$ and $x^\feed$ denote the mole fraction of the most permeable component in the permeate and feed streams, respectively) is constant. One consequence of this assumption is that the permeate composition depends only on the feed composition. This assumption is therefore too strong, because the permeate composition depends significantly on the trans-membrane pressure ratio and the stage cut (fraction of the feed permeated through the membrane).\\ %We remark that \cite{Uppaluri2004, UPPALURI2006832,SCHOLZ20151} consider symmetric membranes, where the bulk permeate composition affects the component fluxes through the membrane. We do not consider such membranes, because they are not prevalent in industrial practice.\\

% To summarize, the crossflow model in conjunction with flux equations derived from the solution-diffusion theory serves as an excellent surrogate for the permeation process. Moreover, it is simple enough to incorporate in a global optimization framework. 

Once the candidate set and a permeator model have been decided, one of the following approaches is used to identify the optimal cascade and its operating conditions (1) explicit enumeration, (2) meta-heuristic approaches such as genetic algorithms or simulated annealing, and (3) mathematical programming approaches. \\

In the first approach, for each cascade, the degrees of freedom are identified, and an exhaustive sensitivity analysis is performed over the admissible operating range of degrees of freedom. The operating condition which optimizes the desired objective is taken to be its optimal operating condition. Next, the cascade which maximizes/minimizes the desired objective at its optimal operating condition is taken to be the optimal cascade (see \cite{BHIDE1993209,BHIDE199113}). The explicit enumeration makes this approach computationally expensive particularly when there are many degrees of freedom for each cascade and the number of cascades is large; for example, cascades containing four or more stages.  \\

% \mt{THIS IS NOT CLEAR:} \jacv{Another limitation of this approach is the following. For a particular set of values to which the degrees of freedom have been set during the sensitivity analysis, if the simulation do not converge to a solution, it is difficult to judge whether the evaluated operating condition is infeasible, or is just that the simulation requires a different initial point to converge. Consequently, after completing the sensitivity analysis, it will not be possible to verify whether the cascade identified as the best one, is truly the global optimal.} \\

In the second approach, simulated annealing or a genetic algorithm is used to determine the optimal cascade and its optimal operating condition (see \cite{Uppaluri2004,UPPALURI2006832,MARRIOTT20034991}). This approach does not provide a global optimality certificate. \\

% Scholz et al. \cite{SCHOLZ20151} use 100 perfect mixers (elements) to obtain a representation of the DAE system.

In the third approach, a mathematical program is formulated and solved using standard solvers (e.g. BARON \cite{tawarmalani2005polyhedral,Kilin2018}, IPOPT \cite{Laird2012}, DICOPT \cite{VISWANATHAN1990769}, etc.). This approach has been used in \cite{QI199871,QI20002719,Aliaga2017,Uppaluri2004,UPPALURI2006832,ADI2016379,KUNDE2018164,SCHOLZ20151}, each of which formulates and solves an MINLP to simultaneously determine the optimal cascade and its optimal operating condition. These formulations use binary variables to model the presence/absence of arcs connecting different membrane stages. In contrast, Pathare and Agrawal \cite{PATHARE2010263} enumerated all cascades in the candidate set explicitly and then formulated a nonlinear program (NLP) for each cascade, which was solved to determine its optimal operating condition. When each cascade in the candidate set had been optimized, the cascade which minimizes the desired objective was taken to be the optimal cascade. The mathematical programming approach suffers from the following challenge. Sensitivity analysis and meta-heuristic approaches can use black-box models, so the first two methods can use higher fidelity models of transport phenomena. On the other hand, mathematical programming approaches typically require explicit functions and/or gradients, and most solvers available today do not allow DAE systems as model constraints. This led prior approaches to resort to the use of discretization procedures to convert the DAE system to a system of algebraic equations. For example, Uppaluri et al.  \cite{Uppaluri2004,UPPALURI2006832} and Scholz et al. \cite{SCHOLZ20151} approximate each membrane stage with several perfect mixers connected in series.  Qi and Henson \cite{QI199871,QI20002719} used a combination of Gauss quadrature and a fourth order Runge-Kutta-Gill method to convert the DAE system into a system of algebraic equations. Kunde and Kienle \cite{KUNDE2018164} discretized each membrane into elements, and in each element, they assumed that the molar flux of each component is linear. This enabled them to compute the permeate and retentate compositions algebraically. All discretization schemes exhibit a trade-off between the model complexity (number of nonlinear equations, etc.) and the accuracy of the solution. The finer the discretization scheme the more accurate is the representation of the DAE system, but the number of nonlinear nonconvex equations is also larger. This makes it challenging to obtain an optimality certificate. Most of the works in the literature use a local solver to solve their formulation. For instance, Qi and Henson \cite{QI199871,QI20002719} use DICOPT++ \cite{VISWANATHAN1990769}, Aliaga et al. \cite{Aliaga2017} use SBB \cite{GAMS2014}, Pathare and Agrawal \cite{PATHARE2010263} use fmincon in MATLAB \cite{MATLAB:2017b}. As we mentioned in \textsection 1, these local solvers can get trapped in suboptimal solutions. We address the above challenges by using the analytical solution of the DAE system, and by using a global solver, such as BARON. We express the analytical solution in a different form than it is known in the literature \cite{weller1950separation,weller1950erratum}, which facilitates the global solution of the optimization problem with off-the-shelf global solvers. We remark that BARON was used in some studies such as those of Adi et al. \cite{ADI2016379}, Scholz et al. \cite{SCHOLZ20151}, and Kunde and Kienle \cite{KUNDE2018164}. However, we note that there are limitations. First, as mentioned before, the perfect mixing model used by Adi et al. \cite{ADI2016379} does not agree well with experimental data, and Scholz et al. \cite{SCHOLZ20151} do not report whether optimality certificate was obtained. The exception is Kunde and Kienle \cite{KUNDE2018164} who do report that they obtain an optimality certificate. Nevertheless, we remark that they treat the trans-membrane pressure difference as a parameter and do not report extensive numerical results with their model. Therefore, the scaling and robustness of their solution procedure is not well understood. \\

% \jacv{As compared to our optimization model, the formulation of Kunde and Kienle \cite{KUNDE2018164} also differs in the objective function. While they minimize the ratio between the sum of permeate flows in the cascade and the feed flow, we chose to optimize the total energy consumption.}\\

% \jacv{Professor Agrawal, you asked what is the method used by Caballero \cite{ARIAS2016371}, and what they miss? Here are the answers. They used the crossflow model for gaseous mixtures and approximate its solution using a finite difference method. However, it does not provide any insight on the accuracy of the approximated model. They use local solvers, and thus, global optimality is not guaranteed. }\\

Finally, Agrawal and Xu \cite{XU1996115,AGRAWAL1996129,XU1996365}, and Pathare and Agrawal \cite{PATHARE2010263} investigated the effect of the exergy loss due to mixing, referred to as \textit{mixing losses}, on the overall compression power. Here, the authors sought the operating condition that minimizes (or completely eliminates in the case of \cite{PATHARE2010263}) the total mixing loss. Empirically, they observed that this approach yields a near optimal solution for cascades containing five or more stages. However, for cascades containing fewer stages, their approach does not always yield the global optimum.  

% \jacv{Empirically, they observed that for a tested example, this approach yielded a near optimal solution for a cascade containing five membrane stages. Such cascade requires four compressors, so it may not be economical. For this reason, industrial practitioners would be interested in cascades containing fewer stages/compressors. In this case, however, their approach does not yield to the global optimum cascade.}

% On the other hand, for cascades containing fewer stages, the optimal operating condition has a substantial mixing loss. \mt{Unclear:} \jacv{As can be noted, there is a balance between the exergy loss due to mixing and the exergy loss within membranes. However, since the relative importance of mixing loss is not known \textit{a priori}, we do not place the restriction of operating with no-mixing losses in this work.}

\section{Permeator model} 
In the rest of the article, we denote the most and the least permeable components as $A$ and $B$, respectively. Further, we denote the separation of a binary mixture $AB$ as $A$/$B$. We use the solution-diffusion theory \cite{WIJMANS19951} to model the local flux of each component through the membrane, and the \textit{crossflow model}, proposed by Weller and Steiner \cite{weller1950separation,weller1950erratum}, to model the overall permeation process.  \\

For both gaseous and liquid mixtures, the solution-diffusion theory is the widely-accepted mechanism of mass transfer through dense polymeric membranes \cite{WIJMANS19951}. According to this theory, the constituent components of the mixture are separated due to their differences in solubility and diffusivity within the membrane. Depending on the type of the mixture, the local flux of each component through the membrane can be obtained using the equations below \cite{WIJMANS19951}.

\begin{subequations}
\begin{align}
& \text{For liquids,} && \left\{
\begin{aligned}
& n_A= \permeance_A \left[x-y\exp\left( -\frac{V_A(P^\ret-P^\per)}{RT}\right)\right]\\
& n_B=\permeance_B \left[(1-x)-(1-y)\exp \left(-\frac{V_B (P^\ret-P^\per)}{RT}\right)\right]
\end{aligned}\right.,
\label{flux_liq}\\
& \text{For gases,} && \left\{\begin{aligned}
    & n_{A} = \permeance_A \left[P^\ret x-P^\per y \right]\\ 
    & n_{B} = \permeance_B \left[P^\ret (1-x)-P^\per (1-y)\right]
    \end{aligned}\right.,
\label{flux_gas}
\end{align}
\end{subequations}

where (i) $n_A$, $\permeance_A$ and $V_A$ (\resp $n_B$, $\permeance_B$, and $V_B$) correspond to the local flux, permeance and molar volume of $A$ (\resp $B$), (ii) $x$ and $P^\ret$ (\resp $y$ and $P^\per$) denote the local mole fraction of $A$ and the total pressure on the retentate (\resp permeate) side, and (iii) $R$ and $T$ denote the universal gas constant and the absolute temperature of the mixture respectively. The expression inside the square brackets in both \eqref{flux_liq} and \eqref{flux_gas} corresponds to the driving force, and it is different for liquids and gases. We derive a unified expression for the local flux of each component that is applicable for both gases and liquids. This enables us to formulate a common optimization model for both gases and liquids. Towards this, we define the following variables:
\begin{subequations}
\begin{align}
    & u=\begin{cases} \ln{r}, &  \text{for a gaseous mixture}\\ 
    \Delta P^{trans}, & \text{for a liquid mixture}\end{cases},
    \label{u_definition}\\
    & \beta = \begin{cases}
    1, & \text{for a gaseous mixture}\\
    0, & \text{for a liquid mixture} \end{cases} ,
    \label{beta-definition}\\
    & C_A=\begin{cases} 1, &  \text{for a gaseous mixture}\\ 
    V_A/RT, &  \text{for a liquid mixture}\end{cases},
    \label{A_definition}\\
    & C_B=\begin{cases} 1, &  \text{for a gaseous mixture}\\ 
    V_B/RT, & \text{for a liquid mixture}\end{cases},
    \label{B_definition}
\end{align}
\end{subequations}
where $r=P^\ret/P^\per$ is the pressure ratio, and $\Delta P^{trans}=P^\ret-P^\per$ is the trans-membrane pressure difference. We now obtain a expression for the local fluxes by expressing \eqref{flux_liq} and \eqref{flux_gas} in terms of $u$, $\beta$, $C_A$ and $C_B$ as 
\begin{subequations}
\label{eq:unified-fluxes}
\begin{align}
& n_A = \permeance_A (P^\ret)^\beta \left[ x- y\, \mathrm{e}^{(-C_A \, u)} \right],\\
& n_B = \permeance_B (P^\ret)^\beta \left[(1-x) - (1-y)\, \mathrm{e}^{(-C_B \, u)} \right].
\end{align}
\end{subequations}

% Here, the permeate travels perpendicularly to the feed flow direction. As a result, the composition in the vicinity of the dense membrane varies only in the feed flow direction on both permeate and retentate sides. When the Assumption 1 is valid, the degree of separation obtained when employing the cross-flow pattern is close to that obtained through an asymmetric membrane where the feed and the bulk permeate flow are either in counter-current (see Figure \ref{Fig_Sing_concep_mem_2}(a)) or in co-current fashion (see Figure \ref{Fig_Sing_concep_mem_2}(b)) \cite{Pan1983}.

\begin{figure}[h]
    \centering
    \includegraphics[width=.5\textwidth]{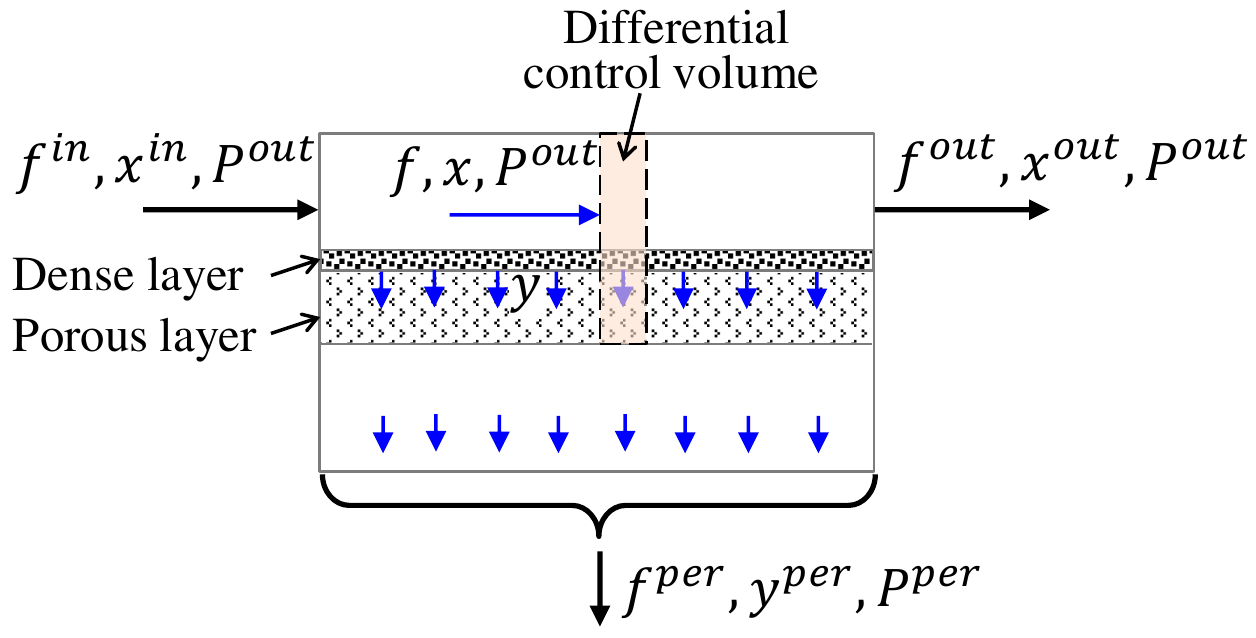}
    \caption{Conceptual representation of crossflow pattern}
    \label{Fig_Sing_concep_mem}
\end{figure}

\begin{figure}[h]
    \centering
    \includegraphics[width=.9\textwidth]{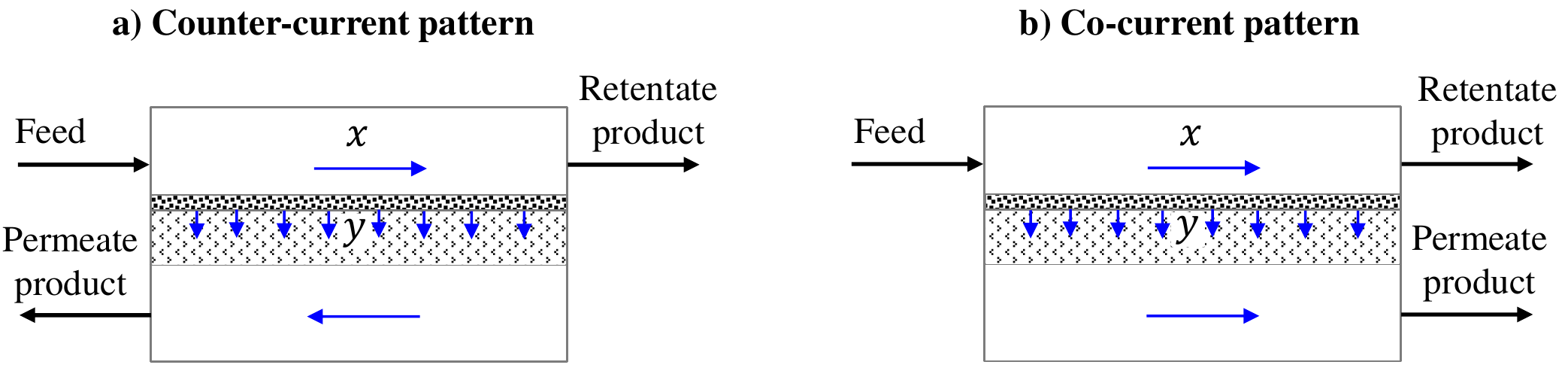}
    \caption{Permeation across an asymmetrical membrane under counter-current and co-current flow patterns. The local composition outside of the dense layer ($y$) does not depend on the flow pattern.}
    \label{Fig_Sing_concep_mem_2}
\end{figure}

Figure \ref{Fig_Sing_concep_mem} shows a schematic of the permeation process across a membrane employing the cross-flow pattern. We make the following assumptions:
\begin{enumerate}
    \item The pressure drop along the membrane module due to the bulk flow of both permeate and retentate streams is not substantial. 
    \item Concentration polarization does not occur near the surface of the membrane.
    \item Mass transfer resistance in the bulk permeate and retentate streams is negligible.
    \item The separation takes place isothermally.
    \item Membrane selectivity is independent of the operating pressure and the composition of the mixture.
\end{enumerate}
Note that, when the Assumption 1 is valid, the degree of separation obtained from an asymmetric membrane employing either the counter-current (see Figure \ref{Fig_Sing_concep_mem_2}(a)) or the co-current (see Figure \ref{Fig_Sing_concep_mem_2}(b)) flow pattern is close to that obtained with the cross-flow pattern \cite{Pan1983}. This is because, the porous layer prevents axial mixing of the local permeate just outside of the dense layer. Consequently, the flux profiles and the net separation remains the same regardless of the flow pattern. Therefore, our models are also applicable for cascades employing either the co-current or the counter-current flow pattern, provided that the Assumption 1 holds and they use asymmetric membranes.\\

%\jacv{Furthermore, for the discussed scenario of small pressure drop at both permeate and retentate sides, the thermodynamic efficiency of the membrane is almost independent of the flow-pattern. Therefore, our single stage membrane model, as well as the membrane cascade optimization model presented in the next section, are also applicable when employing asymmetrical dense membranes operating in co-current or counter-current flow pattern, provided the pressure drop on both permeate and retentate sides is small.} \ra{[Professor Agrawal, let us discuss this part]} 

Let the flowrate, pressure and composition of the mixture entering the membrane module be $f^\feed$, $P^\ret$, and $x^\feed$, respectively. Consider a differential control volume as shown in Figure \ref{Fig_Sing_concep_mem} (dashed rectangle). Mass balance of component $A$ across the control volume is given by $d(fx)=y \,df$, where $x$ (\resp $y$) corresponds to the local mole fraction of component $A$ on the retentate (\resp permeate) side, and $f$ corresponds to the local molar flowrate on the retentate side. The mass balance equation can be rearranged as
\begin{equation} \label{diffpermmodel_1}
\frac{dx}{df} = \frac{y-x}{f},\quad x(f^\feed) = x^\feed,\quad f\in [f^\feed,f^\ret],
\end{equation}

where $x(f^\feed) = x^\feed$ is the initial condition, and $f^\ret$ denotes the net flowrate of the retentate leaving the membrane module. Since there is no axial mixing in the porous layer, the local mole fraction of each component on the permeate side is simply the ratio of the local flux of the component to the total local flux \ie{} $y=n_A/(n_A+n_B)$ and $(1-y)=n_B/(n_A+n_B)$. Since both the equations are linearly dependent, we use only the former. We substitute $n_A$ and $n_B$ from \eqref{eq:unified-fluxes} and rearrange (see \textsection \ref{sec:flux-reformulation} for details) to obtain
\begin{align}
\label{genfluxrel_2}
y-x = k + k (S-1) y - \frac{k\, S}{S-(S-1)y}, 
\end{align}
where $S=\permeance_A/\permeance_B$ is the selectivity of component $A$ w.r.t component $B$ and
\begin{align}
k=\frac{(S-1)-(S e^{-C_A u}-e^{-C_B u})}{(S-1)^2} .
    \label{k_equation}
\end{align}

We solve the DAE system in \eqref{diffpermmodel_1} and \eqref{genfluxrel_2} analytically (see \textsection \ref{sec:analytic-der} for derivation) to obtain
\begin{align}
S\ln \frac{y^\exit}{y^\feed} - \ln \frac{1-y^\exit}{1-y^\feed} - k (S-1)^2\ln \frac{y^\exit - x^\ret}{y^\feed-x^\feed} = k(S-1)^2 \ln(1-\theta),
\label{eq:permeater-model}
\end{align}
where $y^\feed = y(f^\feed)$, $y^\exit=y(f^\ret)$, $x^\ret=x(f^\ret)$, and stage cut $\theta=(f^\feed-f^\ret)/f^\feed$. Mole fractions $y^\feed$ and $y^\exit$ are related to $x^\feed$ and $x^\ret$ via \eqref{genfluxrel_2} \ie{}
\begin{align}
& y^\feed-x^\feed = k+ k(S-1)y^\feed - \frac{k \, S}{S-(S-1)y^\feed}, 
\label{eq:flux-rel-feed}\\
& y^\ret-x^\ret = k+ k(S-1)y^\ret - \frac{k \, S}{S-(S-1)y^\ret}.
\label{eq:flux-rel-ret}
\end{align}
Therefore, given $f^\feed$, $x^\feed$, $r$ (or $\Delta P^{trans}$), and stage cut $\theta$, the mole fraction of component $A$ in the retentate, $x^\ret$, can be determined by solving \eqref{eq:permeater-model}--\eqref{eq:flux-rel-ret} simultaneously. Next, the mole fraction of component $A$ in the permeate, $y^\per$, can be obtained from the overall component mass balance around the membrane module,
\begin{equation}
f^\feed x^\feed=f^\ret x^\ret + f^\per y^\per,
\label{eq:PM-OMB}
\end{equation}
where $f^\per = f^\feed \theta = f^\feed -f^\ret$.

\subsection{Validation of the permeator model}

\begin{figure}[ht]
    \centering
    \includegraphics[width=.8\textwidth]{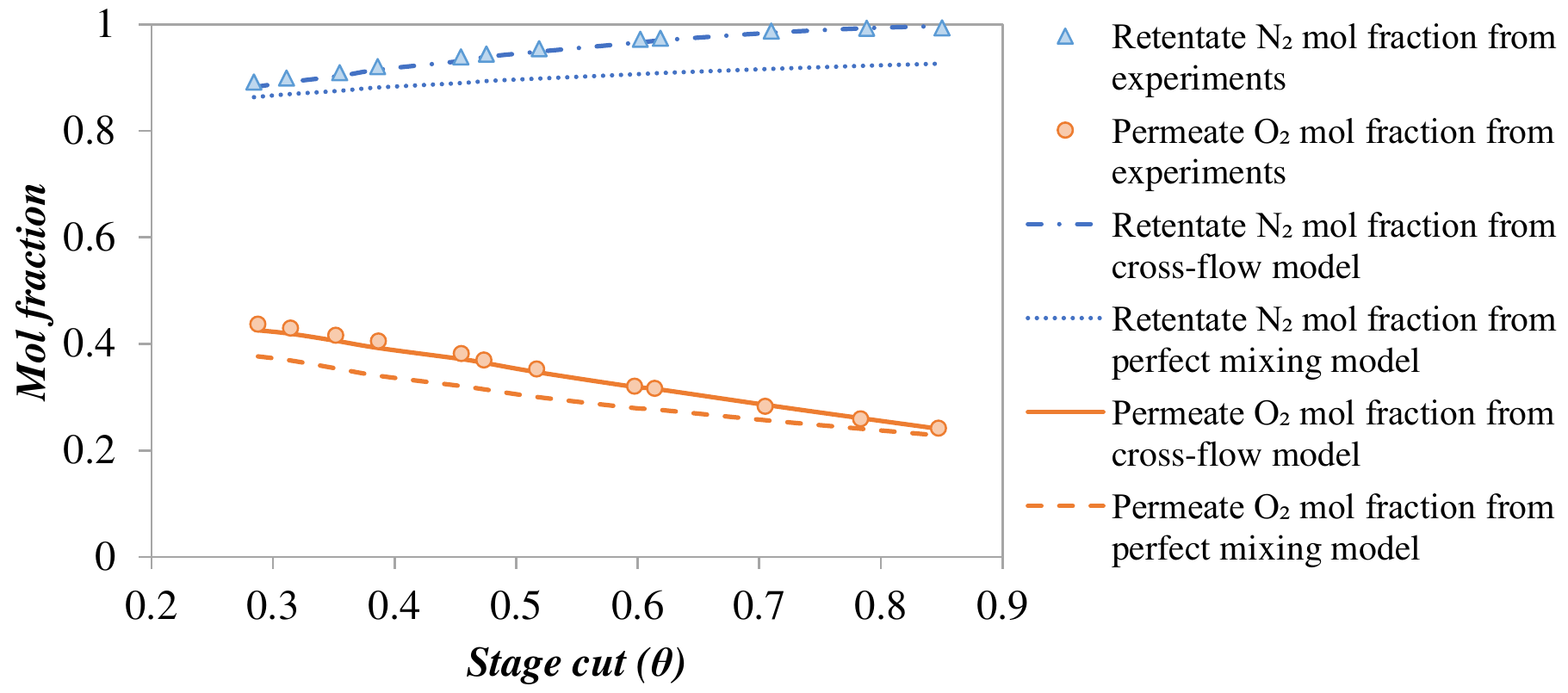}
    \caption{Comparison between experimental data for the separation of \ce{O2}/\ce{N2} \cite{Feng1999} and the predictions obtained from crossflow model and the perfect mixing model. Here, $P^\ret/P^\per = 7.83 \text{ bar}/0.93 \text{ bar} = 8.4$, $x_{\ce{O2}}^\feed=0.205$, membrane permselectivity, $S = 5.3$ (calculated based on the average values of component permeances reported in \cite{Feng1999}).}
    \label{Fig_validation_N2O2}
\end{figure}

\begin{figure}[ht]
    \centering
    \includegraphics[width=.8\textwidth]{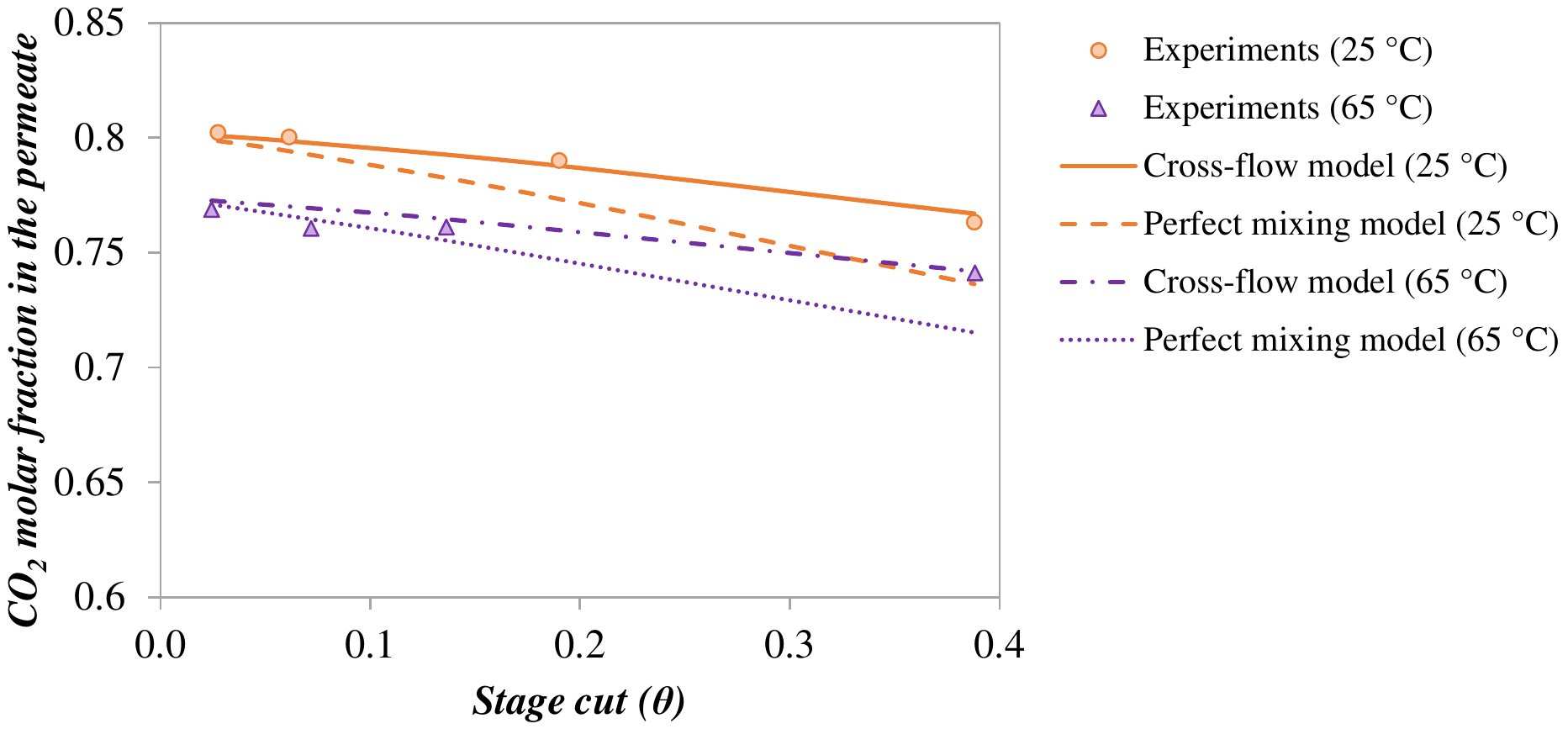}
    \caption{Comparison between experimental data for the separation of \ce{CO2}/\ce{CH4} \cite{Tranchino1989} and the predictions obtained from the cross-flow model and the perfect mixing model. Here, $P^\ret/P^\per = 4.05 \text{ bar}/1.01 \text{ bar} = 4$, $x_{\ce{CO2}}^\feed=0.60$, the membrane permselectivity at 25 \degree C is: 3.58, and the membrane permselectivity at 65 \degree C is: 2.9 \cite{Tranchino1989}.  }
    \label{Fig_validation_CO2CH4}
\end{figure}

We now show the validity of the model by comparing the predicted permeate and retentate mole fractions as a function of stage cut against the experimental data for \ce{O2}/\ce{N2} and \ce{CO2}/\ce{CH4} separations in Figures \ref{Fig_validation_N2O2} and \ref{Fig_validation_CO2CH4}, respectively. In the Figures, we also show the predictions obtained from the \textit{perfect mixing model} \cite{weller1950separation}. Clearly, there is a good agreement between the crossflow model (\eqref{eq:permeater-model}
--\eqref{eq:flux-rel-ret}) and the experimental data for both the mixtures. On the other hand, the perfect mixing model always underestimates the composition of the retentate and the permeate streams leaving the membrane module. Therefore, the optimization results obtained using \eqref{eq:permeater-model}--\eqref{eq:flux-rel-ret} as the permeator model are more reliable than those obtained using the perfect mixing model.\\

We are not aware of literature that provides, for liquid mixtures, detailed experimental data regarding the composition of permeate and retentate streams as a function of stage cut. We will nevertheless use \eqref{eq:permeater-model}--\eqref{eq:flux-rel-ret} for liquid mixtures because the agreement between the local flux determined using \eqref{eq:unified-fluxes} and the experimental value is very good \cite{WHITE2002191,SILVA2010167}. We recognize that this may not sufficiently validate \eqref{eq:permeater-model}--\eqref{eq:flux-rel-ret} for liquid mixtures, but unfortunately the verification of composition profile is not possible in the absence of experimental data.

\pagebreak

\subsection{Minimum Selectivity Requirement}
\label{sec:min-S-require}
By definition, selectivity is the ratio of the permeances of the more permeable and less permeable components, so $S>1$. Now, consider the flux equation in \eqref{genfluxrel_2} in the following form 
\begin{equation}
    y-x = k (S-1)^2 \frac{y(1-y)}{S-(S-1) y}.
\end{equation}
Since $0 \le y \le 1$, the sign of the RHS depends on the sign of $k$. For a gaseous mixture, from \eqref{k_equation}, it can be shown  that $k \ge 0$ regardless of the value of the selectivity, provided $r \ge 1$. This implies $y-x \ge 0$, or the local mole fraction on the permeate side is greater than that on the retentate side. This is consistent with the expected behavior of the permeation process. On the contrary, for a given liquid mixture (\ie{} for a given $C_A$ and $C_B$), $k$ can be negative even for some $S>1$. In that case, $y-x \le 0$, and thus, the permeate from the membrane module is enriched in the less permeable component. This phenomenon, referred to as \textit{negative rejection} \cite{Lonsdale1967,PAUL2004371}, has been observed experimentally in the context of separation of phenol/water via reverse osmosis \cite{Lonsdale1967}. In this work, we avoid negative rejection by choosing
\begin{equation}
    S > \frac{1-e^{-C_B\, u}}{1-e^{-C_A \, u}}.
    \label{eq:min-S1}
\end{equation}
Then, from \eqref{k_equation} and \eqref{genfluxrel_2}, $k$ and $y-x$ are guaranteed to be positive.\\

Further, in a typical permeation process, $y$ increases with increase in $r$/$\Delta P^{trans}$ for a given $x$. Therefore, for a given $x$, we require 
\begin{equation}
    \frac{dy}{dr} = \frac{dy}{du} \cdot \frac{du}{dr} \ge 0
\end{equation}
for a gaseous mixture and 
\begin{equation}
    \frac{dy}{d(\Delta P^{trans})} = \frac{dy}{du} \cdot \frac{du}{d(\Delta P^{trans})} \ge 0
\end{equation}
for a liquid mixture. Since $du/dr \ge 0$ and $du/d(\Delta P^{trans}) \ge 0$ (see \eqref{u_definition}), we require $dy/du \ge 0$. Differentiating \eqref{genfluxrel_2} with respect to $u$ yields
\begin{equation}
    \frac{dy}{du} = \frac{dk}{du} (S-1)^2 \left[\frac{Sx}{y^2} + \frac{1-x}{(1-y)^2}\right]^{-1}.
\end{equation}
Since $0 \le y \le 1$ and $0 \le x \le 1$, $dy/du \ge 0$ only when $dk/du \ge 0$. For a gaseous mixture, from \eqref{k_equation}, it can be verified that $dk/du \ge 0$ regardless of the value of the selectivity. On the other hand, for a liquid mixture, $dk/du \ge 0$ only when
\begin{equation}
    S > \frac{C_B}{C_A} e^{(C_A-C_B)u}.
    \label{eq:min-S2}
\end{equation}
We obtain the minimum selectivity needed by combining \eqref{eq:min-S1} and \eqref{eq:min-S2}, and noting that $u$ is a decision variable, as 
\begin{equation}
S > \max_{u^\lobnd\le u \le u^\upbnd } \left\{\frac{1-e^{-C_B\, u}}{1-e^{-C_A \, u}}, \frac{C_B}{C_A} e^{(C_A-C_B)u}\right\}.
\label{eq:min-Selec}
\end{equation}
Here, $u^\lobnd = (\Delta P^{trans})^\lobnd$ and $u^\upbnd = (\Delta P^{trans})^\upbnd$ (\resp $u^\lobnd=\ln r^\lobnd$ and $u^\upbnd = \ln r^\upbnd$), where $(\Delta P^{trans})^\lobnd$ and $(\Delta P^{trans})^\upbnd$ (\resp $r^\lobnd$ and $r^\upbnd$) denote the lower and upper bounds on the trans-membrane pressure difference (\resp pressure ratio) for the given liquid (\resp gaseous) mixture. Note that, \eqref{eq:min-Selec} reduces to $S > 1$ for a gaseous mixture. On the other hand, the minimum selectivity needed to separate the given liquid mixture depends on the molar volumes of the constituent components ($C_A$ and $C_B$) and on the admissible range of the trans-membrane pressure difference. In this work, we consider only those cases where \eqref{eq:min-Selec} holds. 

\subsection{Properties of the Permeator Model}
\label{sec:properties}
Here, we describe the properties of the permeator model which we will use in \textsection \ref{sec:add-cuts} to derive additional cuts to the MINLP. Provided $\eqref{eq:min-Selec}$ holds,
\begin{itemize}
    \item[P1] The mole fraction of component $A$ on the retentate side decreases along the length of the membrane module. Mathematically, this can be readily shown using \eqref{diffpermmodel_1}. Since $y-x \ge 0$ (see \textsection \ref{sec:min-S-require}) and $f \ge 0$, $dx/df \ge 0$. Further, since $f$ decreases along the length of the membrane module, $x$ also decreases. This implies that $x^\ret \le x^\feed$. 
    \item[P2] The mole fraction of component $A$ in the permeate is at least as high as the mole fraction in the feed. Mathematically, from \eqref{eq:PM-OMB}, $x^\feed$ can be expressed as a convex combination of $x^\ret$ and $y^\per$ \ie \; $x^\feed = (1-\theta)x^\ret + \theta y^\per$. Since $x^\ret \le x^\feed$ from P1, $y^\per \ge x^\feed$.
    \item[P3] The local mole fraction of component $A$ on the permeate side ($y$) increases monotonically with an increase in the local mole fraction on the retentate side ($x$). This is because 
    \begin{align}
        \frac{dy}{dx} = y \left[ x + \frac{k(S-1)^2y^2}{[S-(S-1)y]^2} \right]^{-1}
    \end{align}
    (see \textsection \ref{derivation_dydx} for derivation) is non-negative.
    \item[P4] The local mole fraction of component $A$ on the permeate side ($y$) decreases along the length of the membrane module. This is inferred from properties P1 and P3. Thus, $y^\ret \le y^\feed$. Further, since the overall permeate is the aggregate of permeate from differential membrane elements, its composition $(y^\per)$ satisfies $y^\ret \le y^\per \le y^\feed$.
\end{itemize}

\section{Problem Formulation}

\subsection{Problem statement}
The required input to the problem consists of (i) the molar flowrate and the composition of the feed and the product streams, (ii) efficiencies of compressors and turbines (for gas separations)/pumps and turbocharger (for liquid separations), (iii) membrane perm-selectivity, (iv) the range of admissible operating pressure ratio/trans-membrane pressure difference, (v) temperature of the mixture, and (vi) molar volume of the constituent components for liquid mixtures (see Table \ref{Table_parameters_table}).\\

Given a binary mixture along with all the required inputs, the problem is then to identify the membrane cascade that requires at most $N$ stages and consumes least power for the separation.

\renewcommand{\arraystretch}{1.3}
\begin{table}[!ht]
      \caption{List of input parameters}\label{Table_parameters_table} 
     \begin{center}
        \begin{tabular}{p{4.5cm}p{10cm}}
         \toprule
         \textbf{Symbol}& \textbf{Definition}\\
         \midrule
         $N$ & Maximum number of stages in the cascade\\
         $F$, $F^\per$, $F^\ret$ & Molar flowrate of the given feed mixture, permeate product, and retentate product streams, respectively\\
         $X^F$, $Y^\per$, $X^\ret $ & Mole fraction of component $A$ in the feed mixture, permeate product, and retentate product streams, respectively\\
         $\eta^{comp}$ & Isothermal compressor efficiency \\  
         $\eta^{pump}$ & Pump efficiency \\    
        $\eta^{TC}$ & Turbocharger efficiency \\           
        $V_{A}, V_{B}$ & Liquid molar volume of components $A$ and $B$, respectively (needed only for liquid mixtures)\\ 
        $V^F$ & Liquid molar volume of the feed mixture, calculated as $V^F = X^F V_{A} + (1-X^F) V_{B}$ (needed only for liquid mixtures)\\   
        $V^\ret$ & Liquid molar volume of the retentate product, calculated as $V^\ret = X^\ret V_{A} + (1-X^\ret) V_{B}$ (needed only for liquid mixtures)\\            
        $T$& Absolute temperature of the feed mixture\\        
        $S $ & Membrane permselectivity ($\permeance_A/\permeance_B$ )\\      
        $[ r^\lobnd , r^\upbnd] $ & Admissible range of trans-membrane pressure ratio for a gaseous mixture \\
        $[(\Delta P^{trans})^\lobnd , (\Delta P^{trans} )^\upbnd]$ & Admissible range of trans-membrane pressure difference for a liquid mixture\\    
         \bottomrule           
        \end{tabular}
    \end{center}
\end{table}

% We use superstructures similar to the one in Figure \ref{Fig_superstructure_gases} for all gaseous separations considered in this work. We refer the product stream withdrawn from $\mathcal{P}$ (\resp from $R$) as the \textit{permeate} (\resp \textit{retentate}) \textit{product stream}.

\subsection{Membrane Cascade Superstructure}
\label{sec:superstructure}
\begin{sidewaysfigure}
\includegraphics[width=1\columnwidth]{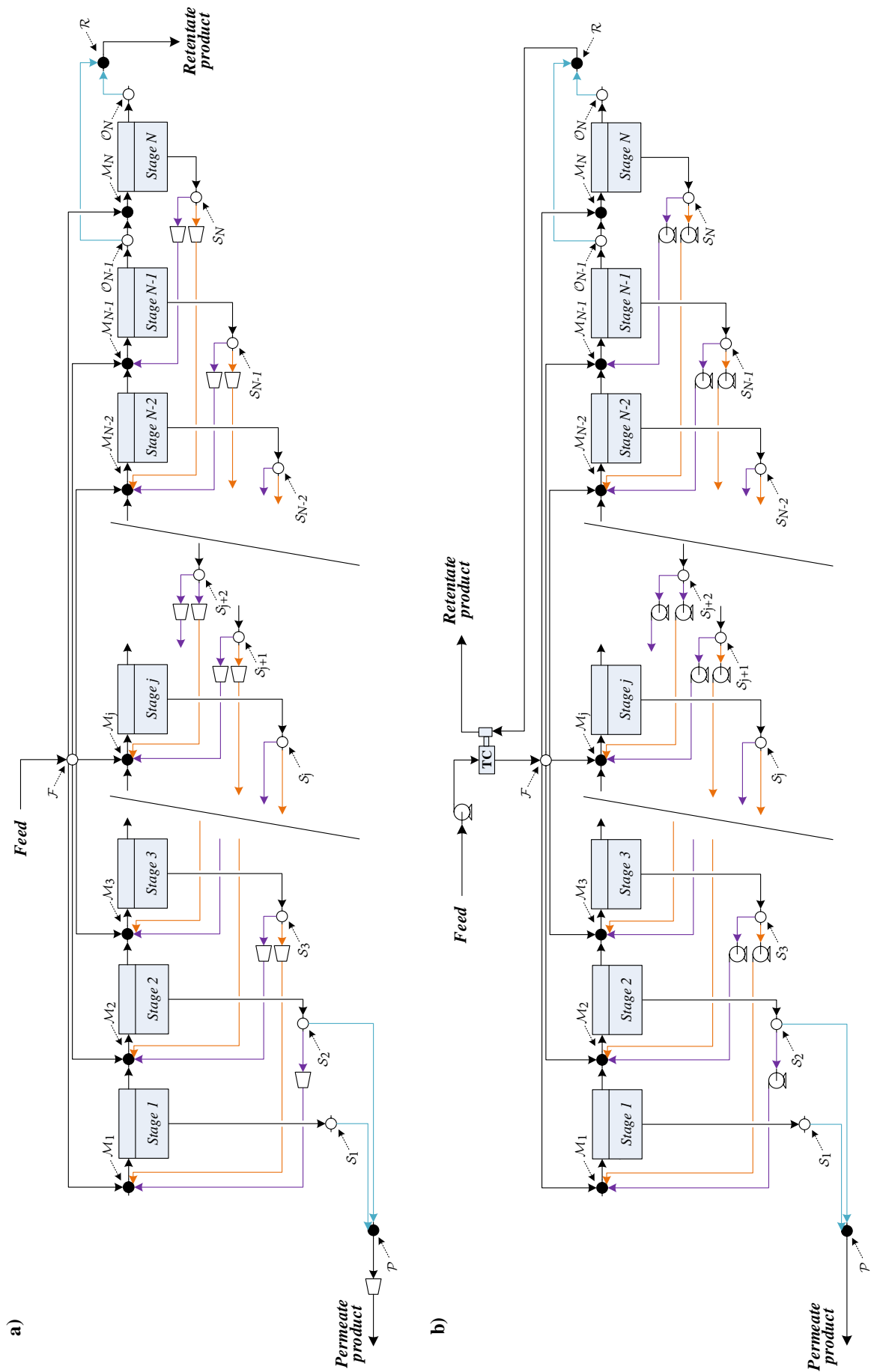}
\caption{Superstructure for the separation of a) gaseous mixtures b) liquid mixtures}
\label{fig:gas-liq-structs}
\end{sidewaysfigure}

Figure \ref{fig:gas-liq-structs}(a) (\resp Figure \ref{fig:gas-liq-structs}(b)) shows the superstructure which embeds cascades requiring at most $N$ stages for the separation of a gaseous mixture (\resp liquid mixture). The splitter $\mathcal{F}$ (see Figure \ref{fig:gas-liq-structs}) splits the feed mixture into $N$ streams which are sent to mixers $\mathcal{M}_1$ through $\mathcal{M}_N$. Each mixer $\mathcal{M}_j$, $j=1,\dots,N-1$, supplies the feed to stage $j$ after mixing the retentate from the stage $j-1$ and the streams from the splitters $\mathcal{F}$, $\mathcal{S}_{j+1}$ and $\mathcal{S}_{j+2}$ (see Figure \ref{fig:gas-liq-structs}). Mixer $\mathcal{M}_N$ supplies the feed to stage $N$ after mixing the streams from the splitters $\mathcal{F}$ and $\mathcal{O}_{N-1}$. The permeate from each stage $j$ is sent to the splitter $\mathcal{S}_j$. Each splitter $S_j$, $j=3,\dots,N$, splits the stream into two streams which are sent to mixers $\mathcal{M}_{j-1}$ and $\mathcal{M}_{j-2}$. Splitter $\mathcal{S}_2$ splits the permeate from stage 2 into two streams which are sent to mixers $\mathcal{M}_1$ and $\mathcal{P}$. Splitter $\mathcal{S}_1$ sends the permeate from stage 1 to the mixer $\mathcal{P}$. On the other hand, the retentate from each stage $j$, $j=1,\dots,N-2$, is sent to the mixer $\mathcal{M}_j$. The retentate from the stage $N-1$ (\resp $N$) is sent to the splitter $\mathcal{O}_{N-1}$ (\resp $\mathcal{O}_N$). Splitter $\mathcal{O}_{N-1}$ splits the inlet stream into two streams which are sent to the mixers $\mathcal{M}_{N}$ and $\mathcal{R}$, and the splitter $\mathcal{O}_N$ sends the retentate from the stage $N$ to the mixer $\mathcal{R}$. Mixer $\mathcal{P}$ (\resp $\mathcal{R}$) mixes the streams from the splitters $\mathcal{S}_1$ and $\mathcal{S}_2$ (\resp $\mathcal{O}_{N-1}$ and $\mathcal{O}_N$) and produces the \textit{permeate product stream} (\resp \textit{retentate product stream}) (see Figure \ref{fig:gas-liq-structs}).\\

In this work, we assume that the given gaseous mixture is at a high pressure and the products are also desired at a high pressure. However, if the gaseous mixture is not available at a high pressure, then an additional compressor can be included in the superstructure to compress the feed. Further, if the products are not desired at a high pressure, then (i) the compressor used for compressing the permeate product stream (see Figure \ref{fig:gas-liq-structs}(a)) can be eliminated from the superstructure, and (ii) a turbine can be included in the superstructure to recover work by expanding the retentate product stream. On the other hand, we assume that the given liquid mixture is always available at a low pressure. We increase the pressure of the feed mixture to an intermediate value using a pump. Next, we send the feed mixture to a turbocharger where it is further pressurized to the desired pressure by transferring the work recovered from the expansion of the retentate product stream  (see Figure \ref{fig:gas-liq-structs}(b)).\\

Further, we impose the following restrictions. 
\begin{enumerate}
    \item All the splitters ($\mathcal{F},\mathcal{S}_1,\dots,\mathcal{S}_N, \mathcal{O}_{N-1},$ and $\mathcal{O}_N$) direct the material flow to one arc entirely. In other words, only one of the arcs arising from a splitter can contain nonzero material flow. 
    \item The trans-membrane pressure ratio/trans-membrane pressure difference is the same for all the stages in the cascade.
\end{enumerate}

Under the above restrictions, we believe that the optimal solution satisfies the following property. 
\begin{itemize}
    \item[P5] In the optimal solution, the mole fraction of the most permeable component in the inlet, retentate and permeate streams decrease from stage 1 through stage $N$ \ie \; $x_1^\feed \ge \dots \ge x_N^\feed$, $x_1^\ret \ge \dots \ge x_N^\ret$, and $y_1^\per \ge \dots \ge y_N^\per$.
\end{itemize}
The above property is based on physical intuition and empirical observation. We use this property in \textsection \ref{sec:add-cuts} to derive additional cuts to the MINLP.

\subsection{Objective function}
\begin{figure}[h!]
    \centering
    \includegraphics[width=0.7\textwidth]{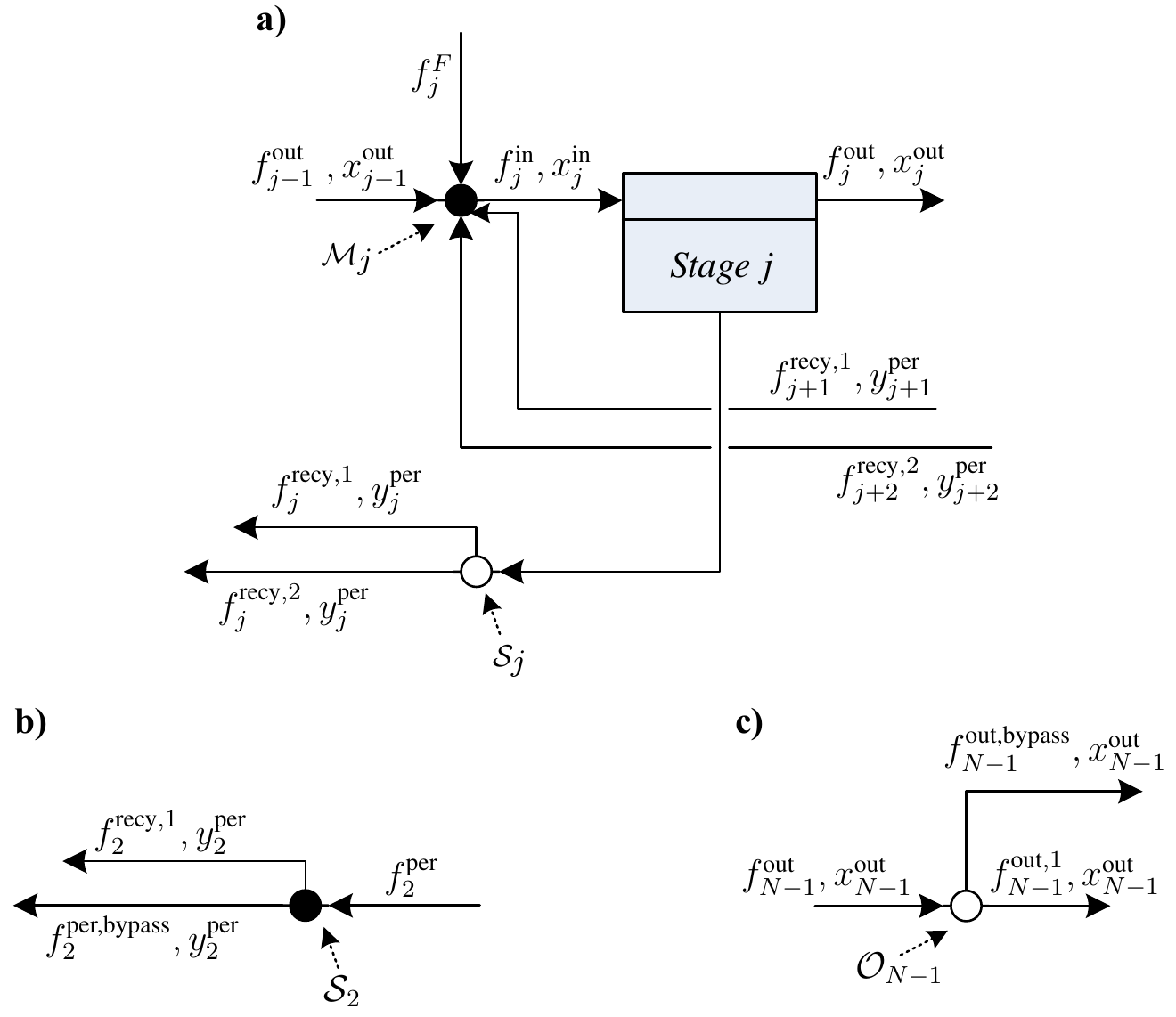}
    \caption{Variables used in the optimization formulation. a) single membrane module, b) splitter $\mathcal{S}_2$, c) splitter $\mathcal{O}_{N-1}$}
    \label{Fig_portion_cascade}
\end{figure}

Figure \ref{Fig_portion_cascade} shows the variables used in our formulation and their definition can be found in Table \ref{Table_variables_table}. 

%Let $N$ denote the maximum number of membrane stages allowed in a cascade. 

\renewcommand{\arraystretch}{1.3}
\begin{table}[!ht]
\caption{List of variables}\label{Table_variables_table} 
     \begin{center}
        \begin{tabular}{p{2cm}p{14cm}}
         \toprule
         \textbf{Symbol}& \textbf{Definition}\\
         \midrule
         $u$& Auxiliary variable defined in \eqref{u_definition} \\     
         $k $ & Auxiliary variable defined in \eqref{perm_const_6}\\    
         $f^F_j$& Molar flow rate along the arc connecting $\mathcal{F}$ and $\mathcal{M}_j$\\
         $f^\feed_j$& Molar flow rate entering stage $j$ from $\mathcal{M}_j$\\
         $f^\ret_j$& Molar flow rate of the retentate stream leaving stage $j$\\
         $f^{\ret,1}_j$& Molar flow rate along the arc connecting $\mathcal{M}_{j+1}$ and $\mathcal{O}_j$\\          
         $f^\bypassret_j$& Molar flow rate along the arc connecting $\mathcal{R}$ and $\mathcal{O}_j$\\         
         $f^\per_j$ & Molar flow rate of the permeate stream leaving stage $j$\\
         $f^{\perrecy,p}_j$ & Molar flow rate along the arc connecting $\mathcal{M}_{j-p}$ and $\mathcal{S}_j$, $p=1,2$ \\
         $f^\bypass_j$ & Molar flow rate along the arc connecting $\mathcal{P}$ and $\mathcal{S}_j$ \\
        $x_j^\feed $ & Mole fraction of component $A$ in the inlet stream to stage $j$\\           
        $x_j^\ret $ & Mole fraction of component $A$ in the retentate stream leaving stage $j$\\ 
        $y^{\per}_j $ & Mole fraction of component $A$ in the permeate stream leaving stage $j$\\
        $y^\feed_j$ & Local mole fraction of component $A$ at the entrance of stage $j$ on permeate side\\         
        $y^\exit_j$ & Local mole fraction of component $A$ at the exit of stage $j$ on permeate side\\
        $z_j^\feed,\; z_j^\exit $ & Auxiliary variables defined in \eqref{opt_permeator}\\      
        $\theta_j $ & Stage cut of stage $j$\\  
        $D_j $ & Auxiliary variable defined in \eqref{eq:Dvar-defs}\\
         $\omega^F_j$ & Binary variable introduced to regulate the flow along the arc connecting $\mathcal{F}$ and $\mathcal{M}_j$ \\ 
         $\omega^{\ret,1}_j$ & Binary variable introduced to regulate the flow along the arc connecting $\mathcal{M}_{j+1}$ and $\mathcal{O}_j$ \\
         $\omega^\bypassret_j$ & Binary variable introduced to regulate the flow along the arc connecting $\mathcal{R}$ and $\mathcal{O}_j$\\  
$\omega^{\perrecy,p}_j$ & Binary variable introduced to regulate the flow along the arc connecting $\mathcal{M}_{j-p}$ and $\mathcal{S}_j$, $p = 1,2 $ \\ 
         $\omega^\bypass_j$ & Binary variable introduced to regulate the flow along the arc connecting $\mathcal{P}$ and $\mathcal{S}_j$\\         
         \bottomrule
        \end{tabular}
    \end{center}
\end{table}

\pagebreak

First, consider the separation of a gaseous mixture. We assume that the given gas mixture is at a high pressure, and that the products are needed at the same pressure as the feed. In this case, energy is needed to compress the outlet streams from splitters $\mathcal{S}_j$ to $r-$times its pressure (see Figure \ref{fig:gas-liq-structs}(a)). We estimate the compressor power using an isothermal compressor efficiency, $\eta^{comp}$. The net power required for the cascade is obtained as
\begin{equation}
    W_G = \frac{R\, T}{\eta^{comp}} \left(\sum_{j=2}^N f_j^{\perrecy,1}+\sum_{j=3}^N f_j^{\perrecy,2}+ F^\per \right) \ln r.
    \label{eq:gas-objective}
\end{equation}

Next, consider the separation of a liquid mixture. We assume that the feed mixture and the product streams are at pressure $P^\per$. In this case, energy is needed to pump the feed mixture from $P^\per$ to an intermediate pressure $P^{TC}$, and to pump the outlet streams from each splitter $\mathcal{S}_j$, $j=2,\dots,N$, from $P^\per$ to $P^\ret$ (see Figure \ref{fig:gas-liq-structs}(b)). Assuming a pump efficiency of $\eta^{pump}$, the net power required for the cascade is obtained as
\begin{align}
W_L = \underbrace{\frac{V^F}{\eta^{pump}} F\, (P^{TC}-P^\per)}_{\text{to pump the feed}} + \underbrace{\frac{1}{\eta^{pump}} \left[\sum_{j=2}^N f_j^{\perrecy,1} V^\per_j\, (P^\ret-P^\per)+ \sum_{j=3}^N f_j^{\perrecy,2} V^\per_j  (P^\ret-P^\per) \right]}_{\text{to pump the recycle streams}}.
\label{eq:liquid-objective}
\end{align}

%\jacv{For the last expression, we have assumed that all streams in the cascade have equal molar volumes, which as described before, were approximated as $V_{mix}=(V_A+V_B)/2$. We note that when the molar volumes of the pure components differ significantly, the previous assumption would not be reasonable, and it would be needed to modify the objective function to account for the molar volume of each stream as function of composition.} 

Here, $V^F$ and $V_j^\per$ denote the molar volume of the feed and the molar volume of the permeate stream from stage $j$, respectively. As mentioned earlier, the turbocharger pressurizes the feed from $P^{TC}$ to $P^\ret$ while expanding the retentate product stream from $P^\ret$ to $P^\per$. Let, $\eta^{TC}$ denote the efficiency of the turbocharger. Then, 
\begin{align}
V^F\, F\, (P^\ret - P^{TC}) = \eta^{TC}\, \left[V^\ret\, F^\ret \, (P^\ret-P^\per)\right].
\label{eq:turbocharger-eff_1}
\end{align}
Here, the LHS corresponds to the power needed to pressurize the feed mixture, the term inside the square brackets on the RHS corresponds to the power produced from the expansion of the retentate product stream, and $V^\ret$ denotes the molar volume of the retentate product stream. We substitute $P^{TC}$ from \eqref{eq:turbocharger-eff_1} in \eqref{eq:liquid-objective} to obtain
\begin{align}
W_L = \frac{1}{\eta^{pump}} \left(F V^F\, \Delta P^{trans} - \eta^{TC} F^\ret V^\ret \, \Delta P^{trans} + \sum_{j=2}^N f_j^{\perrecy,1} V^\per_j \, \Delta P^{trans} +\sum_{j=3}^N f_j^{\perrecy,2} V^\per_j \, \Delta P^{trans} \right).
\label{eq:liq-objec-final}
\end{align}
From \eqref{eq:gas-objective} and \eqref{eq:liq-objec-final}, we obtain a unified objective function that is applicable for both gaseous and liquid mixtures using \eqref{u_definition} as $D_0 u + D_1 F^\per \cdot u + \sum_{j=2}^N D_j f_j^{\perrecy,1} \cdot u + \sum_{j=3}^N D_j f_j^{\perrecy,2} \cdot u + D_{N+1} F^\ret \cdot u$, where
\begin{subequations}
\begin{align}
& D_0 = \begin{cases}
0, & \text{for a gaseous mixture}\\
\frac{V^F\, F}{\eta^{pump}}, & \text{for a liquid mixture}
\end{cases}\\
& D_1 = \begin{cases}
\frac{R\, T}{\eta^{comp}}, & \text{for a gaseous mixture}\\
0, & \text{for a liquid mixture}
\end{cases}\\
& D_j = \begin{cases}
\frac{R\, T}{\eta^{comp}}, & \text{for a gaseous mixture}\\
\frac{V_j^\per}{\eta^{pump}}, & \text{for a liquid mixture}
\end{cases}, \quad j = 2,\dots,N\\
& D_{N+1} = \begin{cases}
0, & \text{for a gaseous mixture}\\
-\frac{V^\ret \, \eta^{TC}}{\eta^{pump}}, & \text{for a liquid mixture}
\end{cases}
\end{align}
\label{eq:Dvar-defs}
\end{subequations}
By assuming ideal mixing, we obtain $V^F = V_A X^F + V_B (1-X^F)$, $V_j^\per = V_A y_j^\per + V_B (1-y_j^\per)$, and $V^\exit = V_A X^\exit + V_B (1-X^\exit)$, where $V_A$ and $V_B$ denote the molar volume of the pure components $A$ and $B$, respectively. Note that $D_0,\dots, D_{N+1}$ are parameters (\resp variables) for a gaseous (\resp liquid) mixture. 

% \jacv{In Equation \eqref{eq:turbocharger-eff_1}, $V_{mix,feed}$ and $V_{mix,retentate}$ represent the molar volume of the feed and retentate product stream respectively. Under the introduced assumption of equal molar volumes for all the streams, we simplify Equation \eqref{eq:turbocharger-eff_1} as shown in Equation \eqref{eq:turbocharger-eff}}

% \begin{align}
% V_{mix}\, F\, (P^\ret - P^{TC}) = \eta^{TC}\, V_{mix}\, f_N^\ret \, (P^\ret-P^\per),
% \label{eq:turbocharger-eff}
% \end{align}

\subsection{MINLP formulation}
Here, we present our mixed-integer nonlinear program (W) for identifying the optimal membrane cascade requiring at most $N$ membrane stages. Consider Figure \ref{Fig_portion_cascade}. Let $\mathcal{J}=\{1,\dots,N\}$. In the following, acronyms OMB and CMB stand for overall mass balance and mass balance of component $A$, $\delta_{(\cdot)} = \{1, \text{ if } (\cdot) \text{ is true}; \; 0, \text{ otherwise}\}$, $(\cdot)^\lobnd$ and $(\cdot)^\upbnd$ denote the lower and upper bounds on $(\cdot)$.
\begin{eqnal}
% Objective function
(\text{W}):\;\;\min \quad & D_0 u + D_1 F^\per \cdot u + \sum_{j=2}^N D_j f_j^{\perrecy,1} \cdot u   \label{Obj_gases_2} \nonumber & & \\
& + \sum_{j=3}^N D_j f_j^{\perrecy,2} \cdot u + D_{N+1} F^\ret \cdot u , & & \text{(Objective Function)}\\
% Definition of D_j
\text{s.t.},\quad & \eqref{eq:Dvar-defs}, & &(\text{Definition of $D_j$})\\
% Overall Mass balance around feed splitter
 & \sum_{j=1}^{N} f_j^{F} =F \label{feed_split_1}, & &(\text{OMB around }\mathcal{F})\\
% Overall mass balance around mixer
\forall \; j \in \mathcal{J},\quad &
\left.\begin{aligned}
f_j^{\feed} = & f_{j-1}^{\ret}\delta_{2 \le j \le N-1} + f_{j-1}^{\ret,1} \delta_{j=N} +f_{j+1}^{\perrecy,1} \delta_{j \le N-1} \\
& +f_{j+2}^{\perrecy,2}\delta_{j \le N-2} +f_j^{F}
\end{aligned}\right\}, & & (\text{OMB around } \mathcal{M}_j)
\label{mass_bal_inlet_1}\\
% Component mass balance around mixer
\forall \; j \in \mathcal{J}, \quad & 
\left.\begin{aligned}
& f_j^{\feed} x_j^{\feed}=  f_{j-1}^{\ret} x_{j-1}^{\ret} \delta_{2 \le j \le N-1} + f_{j-1}^{\ret,1} x_{j-1}^{\ret} \delta_{j=N} \\
& +f_{j+1}^{\perrecy,1} y^{\per}_{j+1} \delta_{j \le N-1} +f_{j+2}^{\perrecy,2} y^{\per}_{j+2} \delta_{j \le N-2}+f_j^{F} X^F
\end{aligned}\right\}, & & (\text{CMB around }\mathcal{M}_j)
\label{mass_bal_inlet_2}\\
% Overall mass balance around retentate splitter
&
\left.\begin{aligned}
& f^{\ret}_{N-1}=f_{N-1 }^{\ret,1}+f_{N-1}^\bypassret\\
& f^\ret_N = f_{N}^\bypassret
\end{aligned}\right\} \label{ret_split_1}, & &(\text{OMB around } \mathcal{O}_j)\\
% Component mass balance around retentate splitter
& 
\left.\begin{aligned}
& f^{\ret}_{N-1} x^\ret_{N-1 }=f_{N-1}^{\ret,1} x^\ret_{N-1} +f_{N-1}^\bypassret x^\ret_{N-1}\\
& f^{\ret}_{N} x^\ret_{N }=f_{N}^\bypassret x^\ret_{N}
\end{aligned} \right\}, & &(\text{CMB around } \mathcal{O}_j)
  \label{ret_split_2}\\
% Overall mass balance around permeate splitter
\forall\; j\in\mathcal{J},\quad & f^{\per}_{j}=f_{j }^{\perrecy,1}\delta_{j \geq 2}+f_{j}^{\perrecy,2}\delta_{j \ge 3}+f_{j}^\bypass \delta_{j \le 2} \label{perm_split_1}, & &(\text{OMB around } \mathcal{S}_j)\\
% Component mass balance around permeate splitter
\forall\; j\in\mathcal{J},\quad &
\left. \begin{aligned}
f^{\per}_{j} y^\per_{j }= & f_{j}^{\perrecy,1} y^\per_{j} \delta_{j \geq 2} + f_{j}^{\perrecy,2} y^\per_{j} \delta_{j \ge 3}\\
& +f_{j}^\bypass y^\per_{j} \delta_{j \le 2}
\end{aligned} \right\}, & &(\text{CMB around } \mathcal{S}_j)
 \label{perm_split_2}\\
% Overall mass balance around stage j
\forall\; j\in\mathcal{J},\quad& f_j^\feed  = f_j^\per +f_j^\ret \label{mass_bal_mem_1}, & &(\text{OMB around stage } j)\\
% Component mass balance around stage j
\forall\; j\in\mathcal{J},\quad & f_j^\feed x_j^\feed = f_j^\per y_j^\per +f_j^\ret x_j^\ret \label{mass_bal_mem_2}, & &(\text{CMB around stage }j)\\
% Overall mass balance around retentate product mixer
& \sum_{j=N-1}^N f_j^\bypassret = F^\ret \label{prod_retmixer_1} & & (\text{OMB around } \mathcal{R})\\
% Component mass balance around retentate product mixer
& \sum_{j=N-1}^N f_j^\bypassret x^\ret_j = F^\ret X^\ret \label{prod_retmixer_2} & & (\text{CMB around } \mathcal{R})\\
% Overall mass balance around permeate product mixer
& \sum_{j=1}^2 f_j^\bypass = F^\per \label{prod_mixer_1} & & (\text{OMB around } \mathcal{P})\\
% Component mass balance around permeate product mixer
& \sum_{j=1}^2 f_j^\bypass y^\per_j = F^\per Y^\per \label{prod_mixer_2} & & (\text{CMB around } \mathcal{P})\\
% Definition of k
& (S-1)^2 k = (S-1)-\left(S e^{-C_A u}-e^{-C_B u}\right), & & (\text{Definition of }k) \label{perm_const_6}\\
\forall \; j \in \mathcal{J},\quad & f_j^\per = \theta_j f_j^\feed, \label{stage_cut} & & \text{(Definition of stage cut)}\\
\forall \; j \in \mathcal{J}, \quad &\left. \begin{aligned}
& S\ln y_j^{\exit}-S\ln y_j^{\feed}-\ln(1-y_j^{\exit})+\ln(1-y_j^{\feed}) \\
& -(S-1)^2(k\ln z_j^{\exit})+(S-1)^2 (k \ln z_j^{\feed})  \\
& =(S-1)^2\left[k\ln(1-\theta_j)\right], \\
& z_j^{\feed}=y_j^{\feed}-x_j^\feed,\\
& z_j^{\exit}=y_j^{\exit}-x_j^\ret,\\
& k+(S-1)(k y_j^{\feed})-S\left[\frac{k}{S-(S-1)y_j^{\feed}}\right]=z_j^{\feed},\\
& k+(S-1)(k y_j^{\exit})-S\left[\frac{k}{S-(S-1)y_j^{\exit}}\right]=z_j^{\exit},
\end{aligned}\right\} & & (\text{Permeator Model})
\label{opt_permeator}\\
& \left. \begin{aligned}
& \omega_j^F & & \in \{0,1\}, & & \forall \; j \in \{1,\dots,N\}\\
& \omega_j^{\perrecy,1} & & \in \{0,1\}, & & \forall \; j \in \{2,\dots,N\}\\
& \omega_j^{\perrecy,2} & & \in \{0,1\}, & & \forall \; j \in \{3,\dots,N\}\\
& \omega_j^\bypass & & \in \{0,1\}, & & \forall\; j \in \{1,2\} \\
& \omega_{N-1}^\bypassret & & \in \{0,1\}\\
& \omega_{N-1}^{\ret,1} & & \in \{0,1\}
\end{aligned} \right\}
\label{eq:bin-definition}\\
% Ensure that the feed is fed to one mixer
& \sum_{j=1}^{N} \omega_j^{F}= 1, 
\label{eq:F-split}\\
\forall \; j \in \mathcal{J}, \quad & \omega^{\perrecy,1}_{j} \delta_{j \geq 2}+\omega^{\perrecy,2}_{j} \delta_{j \ge 3}+\omega^\bypass_{j} \delta_{j \le 2} = 1 
\label{eq:S-split}\\
& \left.\begin{aligned}
& \omega_{N-1}^{\ret,1} + \omega_{N-1}^{\bypassret} = 1,\\
& \omega_N^\bypassret = 1
\end{aligned} \right\},
\label{eq:O-split}\\
& \left. \begin{aligned}
& 0 \le f_j^F \le \omega_j^F (f_j^F)^\upbnd, & & \forall \; j \in \{1,\dots,N\}\\
& 0 \le f_j^{\perrecy,1} \le \omega_j^{\perrecy,1} (f_j^{\perrecy,1})^\upbnd, & & \forall \; j \in \{2,\dots,N\}\\
& 0 \le f_j^{\perrecy,2} \le \omega_j^{\perrecy,2} (f_j^{\perrecy,2})^\upbnd, & & \forall \; j \in \{3,\dots,N\}\\
& 0 \le f_j^{\bypass} \le \omega_j^\bypass (f_j^\bypass)^\upbnd, & & \forall \; j \in \{1,2\}\\
& 0 \le f_j^\bypassret \le \omega_j^\bypassret (f_j^\bypassret)^\upbnd, & & \forall \; j \in \{N-1,N\}\\
& 0 \le f_{N-1}^{\ret,1} \le \omega_{N-1}^{\ret,1} (f_{N-1}^{\ret,1})^\upbnd 
\end{aligned} \right\}
\label{eq:suppress-flows}\\
\forall \; j \in \mathcal{J}, \quad & \left.\begin{aligned}
& 0 \le (\cdot) \le (\cdot)^\upbnd, \quad \forall\; (\cdot) \in \left\{f_j^\feed,f_j^\ret,f_j^\per\right\}\\
& (\cdot)^\lobnd \le (\cdot) \le (\cdot)^\upbnd, \quad \forall \; (\cdot) \in \left\{u, k, \theta_j, x_j^\feed, x_j^\ret,  \right. \\
& \hspace{48mm} \left. y_j^\per, y_j^\feed, y_j^\exit, z_j^\feed, z_j^\exit \right\}
\end{aligned}\right\} & &\text{(Bounds on variables)}\label{variable_bounds}
\end{eqnal}

The formulation of the objective function is described in the previous subsection. We now describe the formulation of constraints.\\

\textbf{Mass balance constraints}: \eqref{feed_split_1} models overall mass balance around the feed splitter $\mathcal{F}$ (see the superstructures in Figure \ref{fig:gas-liq-structs}). \eqref{mass_bal_inlet_1}, \eqref{ret_split_1}, \eqref{perm_split_1}, and \eqref{mass_bal_mem_1} (\resp \eqref{mass_bal_inlet_2}, \eqref{ret_split_2}, \eqref{perm_split_2}, and \eqref{mass_bal_mem_2}) model overall mass balance (\resp mass balance of component $A$) around mixer $\mathcal{M}_j$, splitters $\mathcal{O}_j$ and $\mathcal{S}_j$, and membrane stage $j$, respectively. Mass balances of component $B$ are implied from the overall mass balances and the mass balances on component $A$, so we do not impose them explicitly. However, although \eqref{ret_split_2} (\resp \eqref{perm_split_2}) is implied from \eqref{ret_split_1} (\resp \eqref{perm_split_1}), we impose it explicitly, because it is not implied in the relaxation where the bilinear terms appear in relaxed form. Next, \eqref{prod_retmixer_1} and \eqref{prod_retmixer_2} (\resp \eqref{prod_mixer_1} and \eqref{prod_mixer_2}) model the overall mass balance and the mass balance of component $A$ around the mixer $\mathcal{R}$ (\resp $\mathcal{P}$).  \\

\textbf{Permeator model constraints}: \eqref{perm_const_6} is the same as \eqref{k_equation} and it computes the value of $k$. \eqref{stage_cut} computes the value of the stage cut. \eqref{opt_permeator} is the permeator model in \eqref{eq:permeater-model}--\eqref{eq:flux-rel-ret}. Observe that we introduced auxiliary variables $z_j^\feed$ and $z_j^\exit$ for $y_j^\feed-x_j^\feed$ and $y_j^\exit-x_j^\exit$ along with the constraints $(z_j^\feed)^\lobnd \le z_j^\feed$ and $(z_j^\exit)^\lobnd \le z_j^\exit$ (see \eqref{variable_bounds}), where $0 < (z_j^\feed)^\lobnd$ and $0 < (z_j^\exit)^\lobnd$. The choice of the lower bounds will be discussed shortly. Without the auxiliary variables and the bound constraints, BARON reports an error since it cannot infer $y_j^\feed-x_j^\feed > 0$ and $y_j^\exit-x_j^\exit > 0$, which are needed to well-define $\ln (y_j^\feed-x_j^\feed)$ and $\ln (y_j^\exit-x_j^\exit)$ terms (see \eqref{eq:permeater-model}). Further, we have \textit{disaggregated} all log terms \ie{} each log term of the form $\ln (y_j^\exit/y_j^\feed)$ is expressed as $\ln y_j^\exit - \ln y_j^\feed$. Without disaggregation, a typical factorable relaxation procedure first introduces an auxiliary variable for the fraction $y_j^\exit/y_j^\feed$, and then relaxes the log term over the range of $y_j^\exit/y_j^\feed$. Empirically, we observed that BARON either fails to solve the MINLP, or does so slowly without the disaggregation of the log terms.\\

\textbf{Restricting flows along specific arcs}: As mentioned in \textsection \ref{sec:superstructure}, in this work, we require that all the splitters ($\mathcal{F},\mathcal{S}_1,\dots,\mathcal{S}_N, \mathcal{O}_{N-1},$ and $\mathcal{O}_N$) direct the material flow to one arc entirely. \eqref{eq:bin-definition}--\eqref{eq:suppress-flows} model this requirement. First, we define binary variables (see \eqref{eq:bin-definition} for the domain of index $j$) (i) $\omega_j^F = $ \{1, if the splitter $\mathcal{F}$ directs the material flow to mixer $\mathcal{M}_j$; 0, otherwise\}, (ii) $\omega_j^{\perrecy,1} = $ \{1, if the splitter $\mathcal{S}_j$ directs the material flow to mixer $\mathcal{M}_{j-1}$; 0, otherwise\}, (iii) $\omega_j^{\perrecy,2} = $ \{1, if the splitter $\mathcal{S}_j$ directs the material flow to mixer $\mathcal{M}_{j-2}$; 0, otherwise\}, (iv) $\omega_j^{\bypass} =$ \{1, if the spitter $\mathcal{S}_j$ directs the material flow to mixer $\mathcal{P}$; 0, otherwise\}, (v) $\omega_{j}^{\bypassret} =$ \{1, if the spitter $\mathcal{O}_{j}$ directs the material flow to mixer $\mathcal{R}$; 0, otherwise\}, and (vi) $\omega_{N-1}^{\ret,1} =$ \{1, if the spitter $\mathcal{O}_{N-1}$ directs the material flow to mixer $\mathcal{M}_N$; 0, otherwise\}. Next, \eqref{eq:F-split}--\eqref{eq:O-split} ensure that only one arc, among all the arcs that are leaving from a splitter, is chosen for directing the material flow. \eqref{eq:suppress-flows} suppresses material flow along the arcs when the corresponding binary variable takes the value zero. Here, we choose $(f_j^F)^\upbnd = F$ because the total flowrate along the arcs connecting splitter $\mathcal{F}$ and mixer $\mathcal{M}_j$ cannot exceed the flowrate of the feed. Similarly, since the total flowrate along the arcs connecting $\mathcal{S}_j$ and $\mathcal{P}$ (\resp $\mathcal{O}_j$ and $\mathcal{R}$) cannot exceed the flowrate of the permeate product (\resp retentate product), we choose $(f_j^\bypass)^\upbnd = F^\per$ (\resp $(f_j^\bypassret)^\upbnd = F^\ret$). However, a natural upper bound does not exist on the remaining flow variables. Therefore, we choose a sufficiently large number for $(f_j^{\perrecy,1})^\upbnd, (f_j^{\perrecy,2})^\upbnd, (f_{N-1}^{\ret,1})^\upbnd$. \\

\textbf{Bounds on variables}: It is essential to have finite bounds on all variables, especially those that are involved in nonlinear terms, in order to construct a valid convex relaxation. \eqref{eq:suppress-flows} bounds a few flow variables. \eqref{variable_bounds} bounds the remaining variables in the problem. As before, there is no natural upper bound on $f_j^\feed$, $f_j^\per$ and $f_j^\ret$ variables, so we choose a sufficiently large number for $(f_j^\feed)^\upbnd$, $(f_j^\per)^\upbnd$ and $(f_j^\ret)^\upbnd$. \\

For a gaseous (\resp liquid) mixture, let the admissible range of operating pressure ratio (\resp trans-membrane pressure difference) for the chosen membrane be $[r^\lobnd,r^\upbnd]$ (\resp [$(\Delta P^{trans})^\lobnd,(\Delta P^{trans})^\upbnd$]). Then, we choose the following as the lower and upper bounds on $u$. 
\begin{subequations}
\begin{align}
& u^\lobnd = \begin{cases}
\ln r^\lobnd, & \text{for a gaseous mixture}\\
(\Delta P^{trans})^\lobnd, & \text{for a liquid mixture}
\end{cases}\\
& u^\upbnd = \begin{cases}
\ln r^\upbnd, & \text{for a gaseous mixture}\\
(\Delta P^{trans})^\upbnd, & \text{for a liquid mixture}
\end{cases}
\end{align}
\end{subequations}

We obtain the lower and upper bounds on $k$ by substituting $u=u^\lobnd$ and $u=u^\upbnd$ in \eqref{k_equation}, respectively. Since the chosen perm-selectivity is such that $k$ is a monotonically increasing function of $u$ over the interval $[u^\lobnd,u^\upbnd]$ (see \textsection \ref{sec:min-S-require}), the choice of bounds on $k$ is justified.  \\

By definition, stage cut is the fraction of the total feed permeating through a membrane module, so $\theta_j \in [0,1]$. However, when $\theta_j=1$, $\ln (1-\theta_j)$ (see \eqref{opt_permeator}) is not well-defined. Therefore, we choose $\theta_j^\lobnd = 0$ and $\theta_j^\upbnd = 1-\epsilon_\theta$. For all our computations in this article, we choose $\epsilon_\theta = 10^{-3}$.\\

The choice of upper and lower bounds on mole fraction variables is listed in \eqref{eq:mol-frac-ub-lb}. We choose the composition of the retentate product ($X^\ret$) and the permeate product ($Y^\per$) streams as the lower and upper bound on $x_j^\feed$ (see \eqref{eq:xin-lbub}), respectively. This is justified, because each stage separates the corresponding feed and a further separation would not be needed if the composition of the feed is either above $Y^\per$ or below $X^\ret$. Next, we determine the bounds on $y_j^\feed$ using \eqref{eq:flux-rel-feed} and the bounds on $x_j^\feed$ and $u$. Since $y_j^\feed$ increases monotonically with $x_j^\feed$ (see property P3 in \textsection \ref{sec:properties}) and $u$ (see \textsection \ref{sec:min-S-require}), we obtain $(y_j^\feed)^\lobnd$ (\resp $(y_j^\feed)^\upbnd$) by substituting $x^\feed= (x_j^\feed)^\lobnd$ and $u=u^\lobnd$ (\resp $x^\feed = (x_j^\feed)^\upbnd$ and $u=u^\upbnd$) in \eqref{eq:flux-rel-feed} and solving for $y^\feed$. This procedure is symbolically represented as $y|_{x=(x^\feed)^\lobnd, \; u=u^\lobnd}$ in \eqref{eq:yin-lbub}. Next, from property P5 in \textsection \ref{sec:superstructure}, the mole fraction of component $A$ in the retentate decreases from stage 1 through $N$ \ie \; $x_N^\ret \le \dots \le x_1^\ret$. Further, since the retentate product is formed by mixing the retentate streams from stages $N-1$ and $N$, its composition lies in the interval $x^\ret_{N} \leq X^\ret \leq x^\ret_{N-1}$. Therefore, we impose $X^\ret$ as the lower bound on $x^\ret_{j}$ for every $j\in \{1,\dots,N-1\}$. While zero is a valid lower bound for $x^\ret_N$, specifying it leads to the following issue. When $x^\ret_N =0$, $y^\ret_N=0$ from the flux equation in \eqref{opt_permeator}, and the term $\ln(z^\ret_N)$ in \eqref{opt_permeator} is not well-defined. Therefore, we set $(x^\ret_N)^\lobnd = \epsilon_x$ where $\epsilon_x>0$ (see \eqref{eq:xout-lb}). For all our computations, we choose $\epsilon_x = 10^{-3}$. Next, as before, we obtain $(y_j^\ret)^\lobnd$ using \eqref{eq:flux-rel-ret} and the bounds on $x_j^\ret$ and $u$ (see \eqref{eq:yper-yout-lb}). Next, from property P4 in \textsection \ref{sec:properties}, we have $y_j^\ret \le y_j^\per$. Thus, the lower bound on $y_j^\ret$ is a valid lower bound on $y_j^\per$, so we choose $(y_j^\per)^\lobnd = (y_j^\ret)^\lobnd$ (see \eqref{eq:yper-yout-lb}). Along the same line, the upper bound on $y_j^\per$ is a valid upper bound on $y_j^\ret$, so we choose $(y_j^\ret)^\upbnd = (y_j^\per)^\upbnd$ (see \eqref{eq:yper-yout-ub}). Further, $y_j^\per \le y_j^\feed$ from property P4 in \textsection \ref{sec:properties}, so the upper bound on $y_j^\feed$ is a valid upper bound on $y_j^\per$. However, for every $j\in \{2,\dots,N\} $, a tighter upper bound can be inferred on $y_j^\per$ from the following argument. Since $y_1^\per \ge y_2^\per$ (see property P5 in \textsection \ref{sec:superstructure}), $y_2^\per$ can be at most $Y^\per$ in order to maintain the composition of the permeate product stream at $Y^\per$. Thus, $Y^\per$ is a valid upper bound on $y_2^\per$. Further, $Y^\per$ is also a valid upper bound on $y_3^\per,\dots,y_N^\per$ because $y_2^\per \ge \dots \ge y_N^\per$. This leads to \eqref{eq:yper-yout-ub}. Finally, we obtain the upper bound on $x_j^\ret$ from \eqref{eq:flux-rel-ret} and the bounds on $y_j^\ret$ and $u=u^\lobnd$ as shown in \eqref{eq:xout-ub}. 

\begin{subequations}
\begin{align}
    \forall \; j \in \mathcal{J}, \quad & (x_j^\feed)^\lobnd = X^\ret, \qquad (x_j^\feed)^\upbnd = Y^\per
    \label{eq:xin-lbub}\\
    \forall \; j \in \mathcal{J}, \quad & (y_j^\feed)^\lobnd = y|_{x=(x_j^\feed)^\lobnd, \; u=u^\lobnd}, \qquad (y_j^\feed)^\upbnd = y|_{x=(x_j^\feed)^\upbnd, \; u=u^\upbnd}
    \label{eq:yin-lbub}\\
    \forall \; j \in \mathcal{J}, \quad & (x_j^\ret)^\lobnd = \begin{cases}
    X^\ret, & \text{if } j <N\\
    \epsilon_x, & \text{otherwise}
    \end{cases}
    \label{eq:xout-lb}\\
    \forall \; j \in \mathcal{J}, \quad &  (y_j^\per)^\lobnd = (y_j^\ret)^\lobnd = y|_{(x_j^\ret)^\lobnd, \; u=u^\lobnd}
    \label{eq:yper-yout-lb}\\
    \forall \; j \in \mathcal{J}, \quad  & (y_j^\ret)^\upbnd = (y_j^\per)^\upbnd = \begin{cases}
    (y_j^\feed)^\upbnd, & \text{if }j = 1\\
    Y^\per, & \text{otherwise}
    \end{cases}
    \label{eq:yper-yout-ub}\\
    \forall \; j \in \mathcal{J}, \quad  & (x_j^\ret)^\upbnd = x|_{y=(y_j^\ret)^\upbnd, \; u=u^\lobnd}
    \label{eq:xout-ub}
\end{align}
\label{eq:mol-frac-ub-lb}
\end{subequations}

Lastly, we obtain the bounds on $z^\feed_j$ and $z_j^\exit$ by analyzing the behavior of \eqref{eq:flux-rel-feed} and \eqref{eq:flux-rel-ret}, respectively. It can be verified that in the interval $[0,1]$, the RHS of both the equations is concave, evaluates to zero at $y=0$ and $y=1$, and goes through a maxima at $y= \frac{S-\sqrt{S}}{S-1}$. This leads to the choice of bounds in \eqref{eq:z-bounds}. As before, $z|_{y=(y_j^\feed)^\lobnd,u=u^\lobnd}$ represents that the value of $z$ is obtained by substituting $y=(y_j^\feed)^\lobnd$ and $u=u^\lobnd$ in \eqref{eq:flux-rel-feed}. 
\begin{subequations}
\label{eq:z-bounds}
\begin{align}
    \label{Zin_lb} \forall \; j \in \mathcal{J}, \quad  & (z_j^\feed)^\lobnd=\min \left\{z|_{y=(y_j^\feed)^\lobnd,u=u^\lobnd},\;z|_{y=(y_j^\feed)^\upbnd,u=u^\lobnd} \right\}
    \\
    \label{Zout_lb}
    \forall \; j \in \mathcal{J}, \quad  & (z_j^\exit)^\lobnd=\min \left\{z|_{y=(y_j^\exit)^\lobnd,u=u^\lobnd},\;z|_{y=(y_j^\exit)^\upbnd,u=u^\lobnd} \right\}\\
    \label{Zin_ub}
    \forall \; j \in \mathcal{J}, \quad & (z_j^\feed)^\upbnd=\begin{cases} 
    z|_{y=(y_j^\feed)^\upbnd,u=u^\upbnd}, & \text{if } (y_j^\feed)^\upbnd \le \frac{S-\sqrt{S}}{S-1}\\
    k^\upbnd (\sqrt{S}-1)^2, & \text{if } (y_j^\feed)^\lobnd \leq \frac{S-\sqrt{S}}{S-1} \leq (y_j^\feed)^\upbnd \\ z|_{y=(y_j^\feed)^\lobnd,u=u^\upbnd}, & \text{if }  \frac{S-\sqrt{S}}{S-1} \le (y_j^\feed)^\lobnd\end{cases}\\
    \label{Zout_ub}
    \forall \; j \in \mathcal{J}, \quad & (z_j^\exit)^\upbnd=\begin{cases} 
    z|_{y=(y_j^\exit)^\upbnd,u=u^\upbnd}, & \text{if } (y_j^\exit)^\upbnd \le \frac{S-\sqrt{S}}{S-1}\\
    k^\upbnd (\sqrt{S}-1)^2, & \text{if } (y_j^\exit)^\lobnd \leq \frac{S-\sqrt{S}}{S-1} \leq (y_j^\exit)^\upbnd \\ z|_{y=(y_j^\exit)^\lobnd,u=u^\upbnd}, & \text{if }  \frac{S-\sqrt{S}}{S-1} \le (y_j^\exit)^\lobnd\end{cases}
\end{align}
\end{subequations}

\subsection{Additional cuts}
\label{sec:add-cuts}
Here, we describe the additional constraints that are derived using the properties P1 through P4 in \textsection \ref{sec:properties} and P5 in \textsection \ref{sec:superstructure}. While some of these constraints are redundant to the MINLP (W), they are not implied in the relaxed problem where the nonlinear constraints appear in a relaxed form. Providing these constraints explicitly helps global solvers in expediting the convergence characteristics of the branch-and-bound algorithm. To the MINLP (W), we append
\begin{eqnal}
    \forall\; j\in \mathcal{J},\quad & x^\ret_j\leq x^\feed_j \le y_j^\per,  \label{lin_val_cuts_1}\\
    \forall\; j\in \mathcal{J},\quad & y_j^\exit \le y^\per_j\leq y^\feed_j, \label{lin_val_cuts_2}\\ 
    \forall\; j\in \mathcal{J}\setminus\{1\},\quad & x^\feed_{j}\leq x^\feed_{j-1}, \quad x^\ret_{j}\leq x^\ret_{j-1}, \quad y^\per_j \leq y^\per_{j-1} \label{lin_val_cuts_4}\\
        \forall\; j\in \mathcal{J}\setminus\{1\},\quad & y^\feed_j\leq y^\feed_{j-1}, \quad y^\exit_j\leq y^\exit_{j-1}  \label{lin_val_cuts_5}\\
    & Y^\per \leq y^\per_1 \label{lin_val_cuts_0}\\
    & x^\ret_N \leq X^\ret \label{lin_val_cuts_6}    
\end{eqnal}

We obtain (i) \eqref{lin_val_cuts_1} using the properties P1 and P2, (ii)  \eqref{lin_val_cuts_2} using the property P4, (iii) \eqref{lin_val_cuts_4} using the property P5, and (iv) \eqref{lin_val_cuts_5} using \eqref{lin_val_cuts_4} and the property P3. \eqref{lin_val_cuts_0} and \eqref{lin_val_cuts_6} are derived using the following arguments. Observe that the permeate product stream is a mixture of the permeate streams of stages 1 and 2. Since the upper bound on $y_2^\per$ is $Y^\per$ (see \eqref{eq:yper-yout-ub}), $y_1^\per$ must be at least $Y^\per$ in order to maintain the composition of the permeate product stream at $Y^\per$. Similarly, the retentate product stream is a mixture of the retentate streams of stages $N-1$ and $N$. Since the lower bound on $x_{N-1}^\ret$ is $X^\ret$, $x_N^\ret$ can be at most $X^\ret$ in order to maintain the composition of the retentate product stream at $X^\ret$. This concludes the formulation of the MINLP.

\subsection{Need for a Global Optimization Approach}
\label{sec:global-opt-mot}
\begin{figure}[h!]
    \centering
    \includegraphics[width=.9\textwidth]{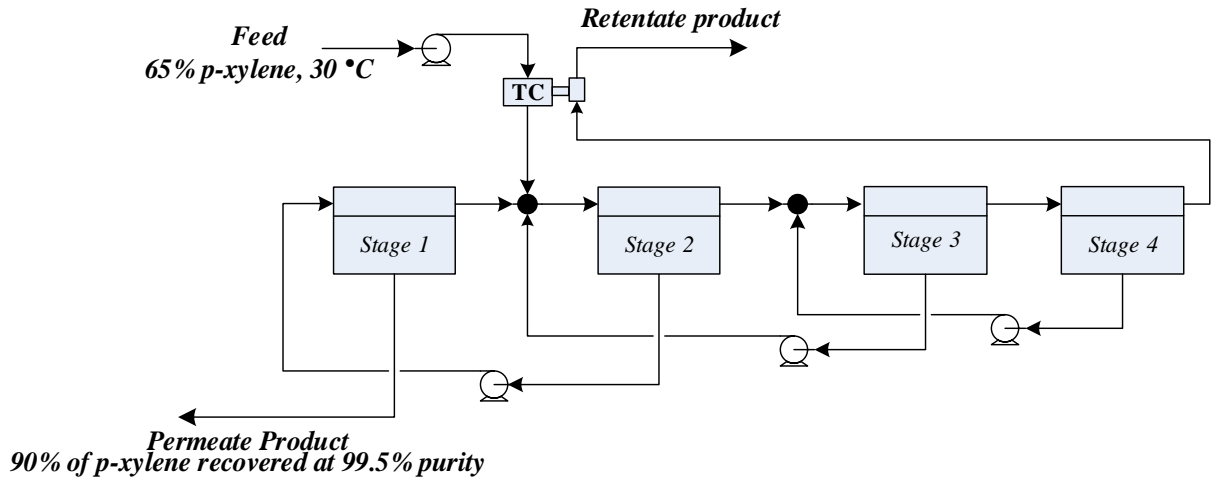}
    \caption{Cascade considered in \textsection \ref{sec:global-opt-mot}. The values of the missing parameters are as follows: $V_A=1.233 \times 10^{-4}$ m\textsuperscript{3}/mol, $V_B=1.215 \times 10^{-4}$ m\textsuperscript{3}/mol, $\eta^{pump}=0.75$, $\eta^{TC}=0.80$, and $[(\Delta P^{trans})^\lobnd , (\Delta P^{trans} )^\upbnd]=$[30 bar, 107 bar]}
    \label{Fig_example_cascade}
\end{figure}

Here, we demonstrate the need for a global optimization approach by showing that local solvers can get trapped in suboptimal solutions even when discrete variables are fixed to a specific cascade. As an example, consider the cascade shown in Figure \ref{Fig_example_cascade} to separate a liquid mixture of p-xylene and a pseudocomponent (mixture of m-xylene and o-xylene). We choose the membrane perm-selectivity to p-xylene to be 50. Let the molar flow rate and the composition (mole fraction of p-xylene) of (i) the feed mixture be 250 mol/s and 0.65, respectively (ii) the permeate product stream be 147 mol/s and 0.995. The molar flowrate and the composition of the retentate product stream can be obtained from an overall mass balance. The values of the remaining parameters are listed in the caption of Figure \ref{Fig_example_cascade}. Our objective is to identify the optimal operating condition of the cascade that minimizes the overall energy consumption. The MINLP (W) can be used for optimization after fixing the binary variables to $\omega_2^F=1$, $\omega^{\perrecy,1}_j = 1$ for $j=2,3,4$, $\omega^{\perrecy,2}_j = 0$ for $j=3,4$, $\omega^\bypass_{2}=0$, and $\omega_{N-1}^\bypassret=0$. When discrete variables are fixed, (W) becomes a nonlinear program, so it can be solved using local solvers such as CONOPT, SNOPT, IPOPT, etc. Here, we use CONOPT. To the best of our knowledge, a systematic method is not available for identifying good initial points. Therefore, we do not provide any initial point to the solver. With its default initialization strategy, CONOPT yields an operating condition that consumes 7,096 kW of power. Whereas, global solver BARON yields an operating condition that consumes 1,780 kW of power ($\sim$75\% reduction in power consumption). This example clearly demonstrates the need to obtain the global optimality certificate when optimizing a membrane cascade. Otherwise, we cannot to be certain whether the obtained solution is the most energy efficient or not.

% \pagebreak

\subsection{Computational experiments}
Here, through numerical experiments, we demonstrate that the proposed MINLP (W) is able to identify the optimal membrane cascade within a relative tolerance of 5\%. In addition, we show the effectiveness of the cuts derived in \textsection 4.5 in expediting the convergence characteristics of BARON by solving (W) with and without \eqref{lin_val_cuts_1}--\eqref{lin_val_cuts_6}. For our numerical experiments, we considered a \textit{test set} of 13 cases described in Table \ref{solverperform}. The values of the remaining parameters are reported in the caption. Note that instead of the flowrate and the composition of permeate and retentate product streams, we have reported the composition ($Y^\per$) and the recovery ($\gamma_A$) of component $A$ in the permeate product stream. Recovery of a component is defined as the ratio of its molar flowrate in the permeate product stream to that in the feed. The values of the parameters needed in \eqref{prod_retmixer_1}-\eqref{prod_mixer_2} can be obtained from
\begin{subequations}
\label{eq:parameters}
\begin{align}
& \gamma_A = \frac{F^\per \cdot Y^\per}{F \cdot X^F},
\label{eq:parameters-1}\\
& F = F^\per + F^\ret,
\label{eq:parameters-2}\\
& F \cdot X^F = F^\per \cdot Y^\per + F^\ret \cdot X^\ret.
\label{eq:parameters-3}
\end{align}
\end{subequations}
\eqref{eq:parameters-1} is the definition of the recovery of component $A$, \eqref{eq:parameters-2} is the overall mass balance across the superstructure, and \eqref{eq:parameters-3} is the mass balance of component $A$ across the superstructure.  For all gaseous (\resp liquid) mixtures, we choose the admissible range of operating pressure ratio (\resp trans-membrane pressure difference) to be $[1.1,9]$ (\resp $[30\text{ bar},107\text{ bar}]$) \ie{} $r^\lobnd = 1.1$ and $r^\upbnd=9$ (\resp $(\Delta P^{trans})^\lobnd = 30$ bar and $(\Delta P^{trans})^\upbnd = 107$ bar).\\

We use BARON 18.5.8 on GAMS 25.1 to solve the MINLP (W). All BARON options except \texttt{pDo} were left at their default values. \texttt{pDo} was set to $-1$. We set the relative tolerance for convergence ($\epsilon_r$) to 5\% and the time limit to 40 h as the termination criteria. All computations are performed on a Dell Optiplex 5040 with 16 GB RAM, which has Intel Core i7-6700 3.4 GHz processor and is running 64-bit Windows 10. The computational results are summarized in Table \ref{solverperform}.\\

The eighth (\resp ninth) column in Table \ref{solverperform} lists the computational performance when the MINLP (W) is solved without (\resp with) \eqref{lin_val_cuts_1}--\eqref{lin_val_cuts_6}. Clearly, when the additional cuts are not included, none of the cases converge even after 40 hours. The remaining duality gap (defined as (Best known upper bound $-$ Best known lower bound)/Best known upper bound) at the end of 40 hours is as high as 98\% in some cases. On the other hand, with the inclusion of \eqref{lin_val_cuts_1}--\eqref{lin_val_cuts_6}, we could solve all 13 cases to 5\%-optimality within 40 hours. It is interesting to note that without \eqref{lin_val_cuts_1}--\eqref{lin_val_cuts_6}, the obtained solutions without optimality certificate, were also at optimality as gleaned by comparing these solutions with those obtained when appending \eqref{lin_val_cuts_1}--\eqref{lin_val_cuts_6}. However, in order to have confidence of weather the obtained solutions are globally optimal, it is essential to obtain the optimality certificate within the desired optimality gap. It avoids cases where one may have suboptimal solutions at hand. Therefore, we recommend to always use constraints \eqref{lin_val_cuts_1}--\eqref{lin_val_cuts_6} when solving the postulated cascade optimization problem.

\begin{sidewaystable}
\caption{Test set for computational experiments. All mole fractions correspond to the most permeable component. The missing parameters are set to $N=4$, $F=250$ mol/s, $\eta^{comp}=0.75$, $\eta^{pump}=0.75$, $\eta^{TC}=0.8$, $V_A=1.233 \times 10^{-4}$ m\textsuperscript{3}/mol, $V_B=1.215 \times 10^{-4}$ m\textsuperscript{3}/mol, and $T = 303.15$ K. Molar volumes are determined using Aspen Plus v8.6, and the reported $V_B$ is the average of molar volumes of o-xylene and m-xylene. All cases were solved with BARON 18.5.8 in GAMS 25.1 with a time limit of 40 hours. Except \texttt{epsr} and \texttt{pDo}, all other BARON options were left at their default values. We set \texttt{epsr} = 0.05 and \texttt{pDo} = -1. TL indicates that the time limit was reached, and the number in parenthesis is the remaining duality gap at the end of 40 hours. I/T correspond to number of BARON iterations/computational time (in hours) to reach 5\%-optimality gap. }
\label{solverperform} 
\centering
\begin{tabular}{ C{1.25cm} C{2.5cm} C{1.5cm} C{2.3cm} C{1.5cm} C{1.5cm} C{1.8cm} C{3cm} C{3cm} }
\toprule
\multirow{3}{*}{Case} & \multirow{3}{*}{Mixture} & Feed & Feed & Permeate & Molar & Membrane & Without & With \\
& & Pressure & Composition & Purity & Recovery & perm- & \eqref{lin_val_cuts_1}--\eqref{lin_val_cuts_6} & \eqref{lin_val_cuts_1}--\eqref{lin_val_cuts_6} \\
& & (bar) & (\% mol) & (\% mol) & (\%) & selectivity & I/T & I/T \\
\midrule
\multicolumn{9}{c}{Gaseous Mixtures}\\
\midrule
1 & \multirow{4}{*}{\ce{CO2}/\ce{CH4}} & 20 & 60 & 95 & 98.7 & 25 & TL (75\%) & 67,523 / 32.8\\
2 & & 20 & 10 & 95 & 81.7 & 25 & TL (89\%) & 14,297 / 3.2\\
3 & & 20 & 10 & 69 & 82.4 & 12 & TL (92\%) & 16,173 / 3.5\\
4 & & 55 & 30 & 86.7 & 93.2 & 21 & TL (74\%) & 12,537 / 3.0\\
\midrule
5 & \ce{H2}/\ce{CO2} & 28 & 21.5 & 99 & 98.2 & 38 & TL (94\%) & 26,627 / 4.8\\
\midrule
6 & \multirow{3}{2cm}{propylene / propane} & 9 & 80 & 90 & 90 & 5 & TL (70\%) & 51,845 / 16.6 \\
7 & & 9 & 70 & 99.6 & 97.8 & 35 & TL (89\%) & 66,903 / 20.8\\
8 & & 9 & 70 & 92 & 97.8 & 35 & TL (62\%) & 33,211 / 13.2\\
\midrule
\multicolumn{9}{c}{Liquid Mixtures}\\
\midrule
9 & \multirow{5}{2cm}{p-xylene / pseudocomponent} & 1 & 65 & \multirow{5}{*}{99.5} & 90 & \multirow{5}{*}{50} & TL (72\%) & 6,631 / 2.2\\
10 & & 1 & 65 & & 99 & & TL (93\%) & 6,963 / 1.8\\
11 & & 1 & 90 & & 90 & & TL (64\%) & 25,848 / 14.5\\
12 & & 1 & 90 & & 99 & & TL (72\%) & 9,091 / 2.7\\
13 & & 1 & 23.6 & & 97.5 & & TL (98\%) & 4,331 / 1.0\\
\bottomrule
\end{tabular}
\end{sidewaystable}

\section{Case studies}
Here, we examine Case 8 (separation of propylene/propane mixture) and Case 12 (separation of p-xylene/(o-xylene+m-xylene) mixture) in Table \ref{solverperform} in more detail. In each case, we determine the optimal membrane cascade requiring at most four stages along with its optimal operating conditions using the MINLP (W). Additionally, we derive constraints that can be appended to the MINLP to limit the number of compressors/pumps in the optimal cascade. All optimization problems in this section are solved to 5\%-optimality gap.

\subsection{Case 8: Propylene/propane separation}

The separation of propylene/propane mixture is of great relevance at industrial level, as propylene is an important precursor in the production of a wide variety of chemicals such as polypropylene, propylene oxide, acrylic acid, acrylonitrile, and isopropanol \cite{propene2013}. Currently, a sub-ambient distillation at mild pressure is used to separate propylene/propane mixtures \cite{ALCANTARAAVILA2014112,FABREGA20101224}. Given the very high demand of propylene, it is typically separated with a high recovery that ranges between 97.8\% to 99\% \cite{colling2004processes,KOOIJMAN2014191}. On the other hand, the final purity of the propylene product varies according to its use in the downstream processes. Propylene needed for the synthesis of polymers is generally produced at $>$99.5\% purity. Whereas, propylene needed for the synthesis of chemicals, such as 2-propanol, is produced at moderate purity of about 92\% \cite{propene2013,propanol2018}.\\

For our case study, we consider the separation of a 70\%/30\% mixture of propylene/propane to produce \textit{chemical grade propylene} at 92\% purity and 97.8\% recovery. The chosen feed composition and product recovery were taken from \cite{colling2004processes}, and they are in good agreement with the values encountered in industry \cite{summers1995high}. We take the molar flow rate of the feed to be 250 mol/s, and we assume that the feed is available at 9 bar. Further, we assume that all the membranes in the cascade have a perm-selectivity of 35 \cite{colling2004processes} (note that propylene is the most permeable component).

% In addition, we assume the molar flow rate of the feed mixture to be 250 mol/s, and that it is available at 9 bar. Lastly, we assume that all the membranes in the cascade have a perm-selectivity of 35 \cite{colling2004processes} (note that propylene is the most permeable component). 

\subsubsection{Optimal four-stage cascade}
\begin{figure}[h!]
    \centering
    \includegraphics[width=1\textwidth]{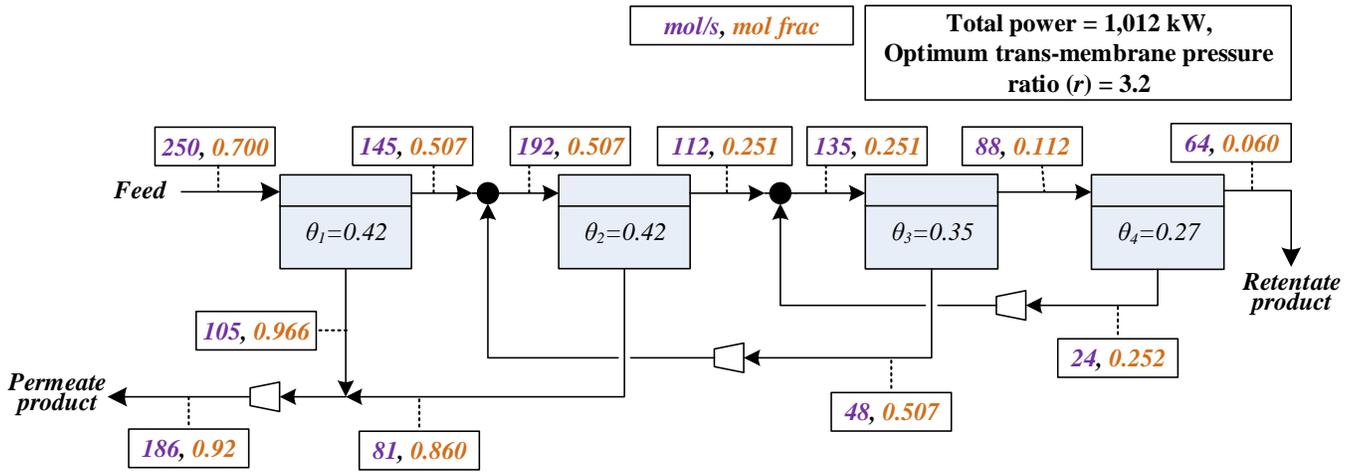}
    \caption{Optimal four-stage cascade for the separation of the propylene/propane mixture. All mole fractions shown in the box-legends correspond to propylene. Any discrepancy in the mass balances is due to the rounding of the flow and the composition values.}
    \label{Fig_cascades_propylene_propane_4stg}
\end{figure}

The optimal four-stage membrane cascade for the separation of the propylene/propane mixture is shown in Figure \ref{Fig_cascades_propylene_propane_4stg}. In this cascade, the feed mixture is located at the first stage. This is not surprising given that the degree of enrichment of the most permeable component in the permeate product stream is not too high. Thus, a single stage is sufficient to achieve the desired purity specification. Nevertheless, it is not sufficient to recover the desired amount of propylene while maintaining the purity. As shown in Figure \ref{Fig_cascades_propylene_propane_4stg}, the optimal cascade uses three additional stages for stripping propylene from the retentate stream leaving the first stage. At optimal operation, the purity of the permeate from the first stage is higher than the desired purity (see Figure \ref{Fig_cascades_propylene_propane_4stg}). Thus, the permeate from the second stage, which does not satisfy the purity requirement (see Figure \ref{Fig_cascades_propylene_propane_4stg}), can be mixed with that from the first stage. This enables the production of the permeate product stream at the desired purity while avoiding the recycling of the permeate from the second stage. 

%The optimum four-stage membrane cascade for the separation of the propylene/propane mixture is shown in Figure \ref{Fig_cascades_propylene_propane_4stg}. In this cascade, the feed mixture is located at the first stage. This is not surprising given that the degree of enrichment of the most permeable component in the permeate product stream is not too high. Thus, one stage is sufficient to achieve the required purity specification. It is also because of this reason, the optimal cascade does not recycle the permeate stream of stage 2 to stage 1. This helps in reducing the energy requirement as it avoids compressing extra material flow. Clearly, in this case, a single membrane stage is sufficient to produce propylene at the desired purity. Nevertheless, a single stage is not sufficient to recover the desired amount of propylene while maintaining the purity. As shown in Figure \ref{Fig_cascades_propylene_propane_4stg}, the optimal arrangement employs three additional stages for stripping propylene from the retentate stream leaving the first stage. \ra{Discuss comment from Prof Agrawal}

\subsubsection{Optimal cascade with only one intermediate compressor}

The optimal four-stage cascade has two intermediate compressors, and one compressor to pressurize the permeate product stream. Owing to the high capital cost and maintenance issues associated with compressors, process designers often seek cascades with few compressors. In the following, we formulate the constraints which can be appended to the MINLP (W) to obtain cascades with a specified maximum number of compressors.  \\

Excluding the compressor needed to compress the permeate product stream, the maximum number of compressors, $M_{comp}$, in a cascade can be controlled using the binary variables $\omega_j^{\perrecy,1}$ and $\omega_j^{\perrecy,2}$ through 
\begin{align}
    \sum_{j=2}^N \omega_j^{\perrecy,1} + \sum_{j=3}^N \omega_j^{\perrecy,2} - \sum_{j=3}^N \omega_{j-1}^{\perrecy,1}\omega_j^{\perrecy,2} \leq M_{comp},
    \label{const_comp0}
\end{align}
The LHS computes the number of compressors in the cascade. The first (\resp second) term in the above equation accounts the compressor along the arc connecting the splitter $\mathcal{S}_j$ and mixer $\mathcal{M}_{j-1}$ (\resp $\mathcal{M}_{j-2}$) (see Figure \ref{fig:gas-liq-structs}(a)). In cases where the material flow from both $\mathcal{S}_{j-1}$ and $\mathcal{S}_j$ is sent to mixer $\mathcal{M}_{j-2}$, both the streams can be mixed before compressing. This reduces the number of compressors needed for the cascade. Thus, the third term in the above equation reduces the number of compressors by one, when both $\omega_{j-1}^{\perrecy,1} = 1$ and $\omega_j^{\perrecy,2} = 1$. We linearize \eqref{const_comp0} by introducing an auxiliary variable $h_{j-1}$ for the product, and relax $h_{j-1}=\omega_{j-1}^{\perrecy,1} \omega_j^{\perrecy,2}$ by replacing the bilinear term with its convex and concave envelopes over $[0,1]^2$. This results in the following constraints
\begin{eqnal}
   & \sum_{j=1}^N \omega_j^{\perrecy,1} + \sum_{j=3}^N \omega_j^{\perrecy,2}-\sum_{j=3}^{N} h_{j-1}\leq M_{comp},  \label{const_comp1}\\
    & \forall\; j\in \left\{ 3,\dots,N\right\},\quad \left\{ 
    \begin{aligned}
    & h_{j-1} \le \omega_{j-1}^{\perrecy,1}, \quad h_{j-1} \le \omega_{j}^{\perrecy,2}\\
    & h_{j-1} \ge 0, \quad h_{j-1} \ge \omega_{j-1}^{\perrecy,1}+\omega_{j}^{\perrecy,2}-1,
    \end{aligned}
    \right.
    \label{const_comp2}
\end{eqnal}
Since $\omega_{j-1}^{\perrecy,1}$ and $\omega_j^{\perrecy,2}$ are binary, it can be easily verified that the feasible solutions to \eqref{const_comp1} and \eqref{const_comp2} are the same as that of \eqref{const_comp0}.\\

Suppose that a stage $j$ does not carryout any separation \ie its stage cut is zero ($\theta_j=0$). Then, the molar flow rate of the permeate stream from stage $j$ is zero, so the compressors along the arcs connecting the splitter $\mathcal{S}_j$ and the mixers $\mathcal{M}_{j-1}$ and $\mathcal{M}_{j-2}$ are not needed. However, \eqref{eq:S-split} does not permit the optimizer to set both $\omega_j^{\perrecy,1}=0$ and $\omega_j^{\perrecy,2}=0$. Instead it requires that either $\omega_j^{\perrecy,1}=1$ or $\omega_j^{\perrecy,2}=1$, which contributes towards the number of compressors in the cascade. Therefore, we modify \eqref{eq:S-split} as shown in \eqref{const_comp3}. Here, $\mu_j$ is a binary variable which takes the value one when the stage cut of a stage is nonzero (see \eqref{const_comp4}). This way, the optimizer can suppress the compressors when the corresponding stage cut becomes zero by setting $\mu_j=0$. 

%Observe that because of constraints \eqref{perm_split_5} and \eqref{perm_split_6}, for each membrane stage beyond the second stage, either $\omega_j^{\perrecy,1}$ or $\omega_j^{\perrecy,2}$ are equal to one, regardless of whether the permeate flow is zero. This implies in the first two summation terms of \eqref{const_comp1} to count one compressor for each membrane even when its permeate flow is zero, which could make the MINLP to be infeasible when adding \eqref{const_comp1}-\eqref{const_comp2} and specifying a relative low value of maximum number compressors. To resolve this issue, we substitute \eqref{perm_split_5} by \eqref{const_comp3}-\eqref{const_comp4}, which allows to set both $\omega_j^{\perrecy,1}$ and $\omega_j^{\perrecy,2}$ to zero when $\theta_j$ is zero, or equivalently, when the permeate flow from stage $j$ is zero. 

\begin{eqnal}
\forall\; j\in\mathcal{J},\quad & \left. \begin{aligned}
& \omega^{\perrecy,1}_{j} \delta_{j \geq 2}+\omega^{\perrecy,2}_{j} \delta_{j \ge 3}+\omega^\bypass_{j} \delta_{j \le 2} = \mu_j,\\
& \omega_j^{\perrecy,1} \in \{0,1\}, \quad \omega_j^{\perrecy,2} \in \{0,1\}, \quad \omega_j^{\bypass} \in \{0,1\} , \quad \mu_j \in \{0,1\} 
\end{aligned} \right\} & &
\label{const_comp3}\\
\forall\; j\in\mathcal{J},\quad & \theta_j \le \mu_j  \label{const_comp4}, & & 
\end{eqnal}

\begin{figure}[h!]
    \centering
    \includegraphics[width=0.9\textwidth]{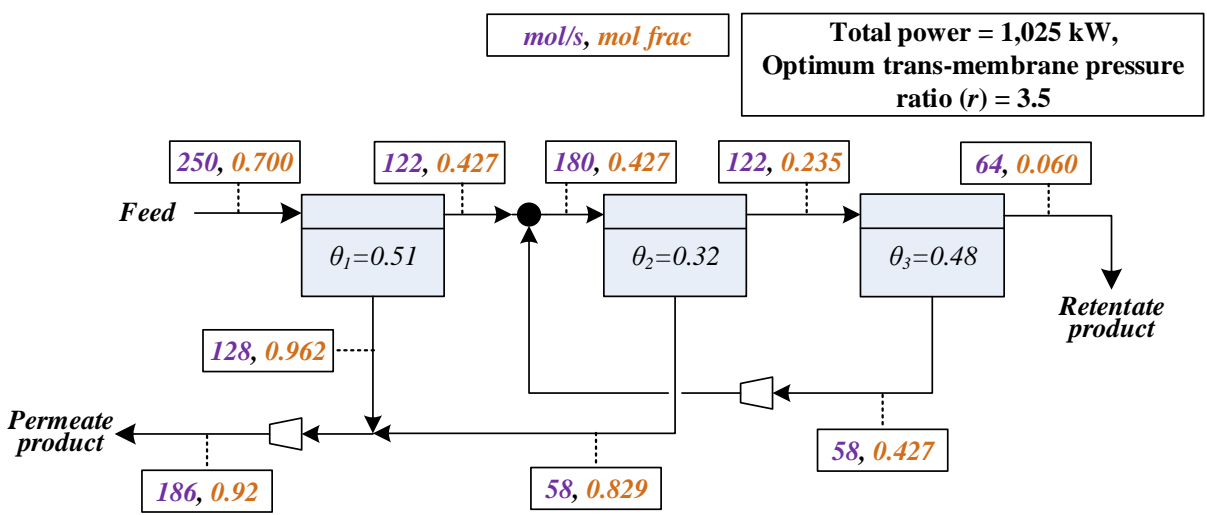}
    \caption{Optimal cascade with only one intermediate compressor for the separation of the propylene/propane mixture. All mole fractions shown in the box-legends correspond to propylene. Any discrepancy in the mass balances is due to the rounding of the flow and the composition values.}
    \label{Fig_cascades_propylene_propane_3stg}
\end{figure}

We now append \eqref{const_comp1}-\eqref{const_comp4} to (W), eliminate \eqref{eq:S-split}, and solve the resulting MINLP by substituting $M_{comp} = 1$ to obtain the optimal cascade requiring only one intermediate compressor. The optimal solution, shown in Figure \ref{Fig_cascades_propylene_propane_3stg}, requires only three stages. The optimizer bypasses material flow through a stage by setting its stage cut to zero. This cascade requires only 1.3\% more power than the one in Figure \ref{Fig_cascades_propylene_propane_4stg}. Since the capital cost of a compressor can be significant, the small increase in energy cost will be dominated by the reduction in the overall capital cost making this arrangement more economical compared to the one in Figure \ref{Fig_cascades_propylene_propane_4stg}. A quantitative cost benefit analysis of the cascades in Figures \ref{Fig_cascades_propylene_propane_4stg} and  \ref{Fig_cascades_propylene_propane_3stg} is beyond the scope of this paper.  \\

% The main purpose of this example is to illustrate the application of constraints \eqref{const_comp1} and \eqref{const_comp2} in identifying attractive cascades that require fewer compressors. 

% We remark that the solution in Figure \ref{Fig_cascades_propylene_propane_3stg} can also be obtained by solving the MINLP (W) with $N=3$, instead of solving the MINLP by appending \eqref{const_comp1} and \eqref{const_comp2}. However, this is a coincidence. As we shall demonstrate with the next case study, in general, solving the MINLP with \eqref{const_comp1} and \eqref{const_comp2} yields better solutions than simply solving the MINLP with fewer number of stages. \\

\subsection{Case 12: Separation of a mixture of xylenes}
Xylenes are important compounds used in the manufacturing of a variety of products ranging from plastics, paint solvents, resins, inks, and plasticizers \cite{daramola2013xylenes}. Among all the isomers, p-xylene has the highest commercial demand. One of the major applications of p-xylene is the synthesis of terephthalic acid, which is an important building block for manufacturing polymers such as polyethylene terephthalate (PET) and polybutylene terephthalate (PBT) \cite{PEREGO201097}. p-Xylene is primarily produced by catalytic reforming of naphtha. In addition, it is also produced via the toluene disproportionation process. However, compared to catalytic reforming, the latter process is more attractive from the following perspective. The C8 cut obtained from the reforming process contains only 20 mol\% of p-xylene. The rest of the mixture consists of o-xylene, m-xylene and ethyl benzene (EB). On the other hand, the C8 cut obtained from the toluene disproportionation process contains around 90 mol\% of p-xylene\cite{ashraf2013analysis,TSAI1999355}, and EB is practically absent in the cut \cite{daramola2013xylenes}. In either case, a separation is needed to concentrate p-xylene since it is needed at $> 99\%$ purity for the downstream processes. \\

Currently, simulated moving bed and crystallization are predominantly used to separate p-xylene from the C8 cut. The inherent energy intensive nature of the separation is due to similarities in physical and structural properties of the constituent components, and it motivates the development of alternative energy-efficient technologies. Recently, reverse osmosis is proposed as a plausible alternative \cite{koh2016reverse}, so we investigate the separation of p-xylene from the C8 cut using a membrane cascade. We assume that the feed mixture comes from the toluene disproportionation process, so the composition of ethyl benzene in the mixture is negligible. Further, we treat the C8 cut as a binary mixture by lumping both o-xylene and m-xylene as a single pseudo component. This is justified, because the membranes available in the literature are more selective to p-xylene compared to both o-xylene and m-xylene \cite{koh2016reverse,CHAFIN20103462}. \\

For our case study, we consider a 250 mol/s mixture of xylenes containing 90\% p-xylene at 30 \degree C and 1 bar. Since p-xylene is a valuable chemical and it is needed at a high purity for downstream processes, we require that the permeate product stream recovers 99\% of p-xylene at 99.5\% purity \cite{TSAI1999355}. The perm-selectivity of p-xylene w.r.t o-xylene+m-xylene pseudo component is taken to be 50. This value is close to the experimentally observed selectivity of p-xylene w.r.t o-xylene through a carbon molecular sieve membrane (see Koh et al. \cite{koh2016reverse}). Finally, we let the admissible operating range of the trans-membrane pressure difference be 30 bar -- 107 bar.

\subsubsection{Optimal four-stage cascade}

\begin{figure}[h!]
    \centering
    \includegraphics[width=0.95\textwidth]{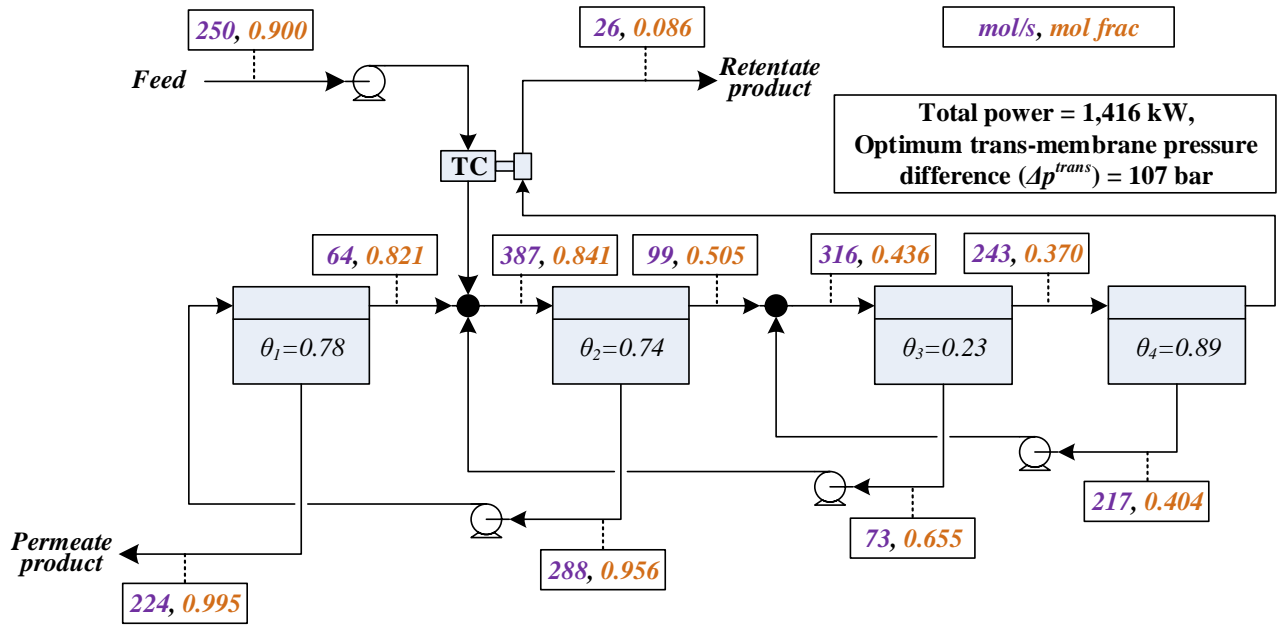}
    \caption{Optimal four-stage cascade for the separation of xylene isomers. All mole fractions shown in the box-legends correspond to p-xylene. Any discrepancy in the mass balances is due to the rounding of the flow and the composition values.}
    \label{Fig_cascades_px_ox_4sgt}
\end{figure}

Figure \ref{Fig_cascades_px_ox_4sgt} shows the optimal four-stage cascade along with its optimal operating condition. In this cascade, the feed mixture is fed to the second stage. With the exception of the first stage, the permeate stream from each stage is pumped and recycled to the stage on its left. For this reason, the cascade requires three pumps to pressurize each of the recycle streams. This is in addition to the pump and the turbocharger needed to pressurize the feed mixture (see Figure \ref{Fig_cascades_px_ox_4sgt}). \\

\subsubsection{On the admissible operating range of $\Delta P^{trans}$}

In Figure \ref{Fig_cascades_px_ox_4sgt}, at optimal operation, we observe that the trans-membrane pressure difference ($\Delta P^{trans}$) is at its upper bound. This suggests that the overall energy consumption could be reduced by increasing the upper bound on $\Delta P^{trans}$. To verify this hypothesis, we increase $u^\upbnd$ and solve the MINLP. The results are shown in Table \ref{Table_sensibilities_px_ox}. As expected, the total power requirement reduces when the cascade is allowed to operate at a higher $\Delta P^{trans}$. Unfortunately, the membranes that are currently available suffer from compaction and lose selectivity towards p-xylene dramatically approximately beyond 107 bar\footnote{In the experimental permeation data provided by Prof. Ryan Lively, we observed that membrane perm-selectivity drops when the trans-membrane pressure difference exceeds 107 bar.}\cite{koh2016reverse}. Therefore, novel membrane materials with improved mechanical properties are needed to further reduce the power requirement. 

% As $\Delta P^{trans}$ is increased, the molar flowrate of the recycle streams decreases leading to the reduction in the overall pump power. 

\begin{table}[!ht]
     \centering
        \begin{tabular}{C{3cm} C{4.5cm} C{3cm}}
         \toprule
         $u^\upbnd$ (bar) &  Optimum $\Delta P^{trans}$ (bar) & Total power (kW)\\
         \midrule
         107   & 107 & 1,416\\
         125   & 125 & 1,281\\
         150   & 150 & 1,205\\
         \bottomrule
        \end{tabular}
        \caption{Optimal power consumption for the cascade in Figure \ref{Fig_cascades_px_ox_4sgt} for different ranges of admissible values of $\Delta p^{trans}$.}
        \label{Table_sensibilities_px_ox}  
\end{table}

\subsubsection{Determining target membrane selectivity}

An alternative, and well-known, approach to reduce the energy requirement of a cascade is to increase the perm-selectivity of the membrane. Currently, the membrane literature relies heavily on the Robeson plot to assess the performance of a membrane \cite{ROBESON2008390}. Unfortunately, this analysis does not answer how the energy requirement of a membrane cascade compares with an alternate technology. Especially, if membranes are being developed in order to replace an existing technology, it is useful to determine the minimum perm-selectivity needed to compete with the existing alternate technology in terms of the energy consumption. We now show how to determine the target selectivity for a given application using our MINLP formulation. \\

We solve the MINLP (W) for different values of selectivity and obtain the optimal power requirement for each value. The results are summarized in Table \ref{Table_selectivities_px_ox}. For the separation of a mixture of xylenes, the state-of-the-art polymeric membranes available in the literature have a selectivity of $\sim$50. For such membranes, the optimal power requirement is 1,416 kW. If the existing non-membrane separation technology consumes more power, then replacing it with a membrane cascade would be beneficial from the point of view of energy consumption. On the other hand, if the existing separation technology consumes less power, say by 10\%, then replacing it with a membrane cascade would be detrimental. In that case, a membrane perm-selectivity of at least 63 is needed to compete with the existing technology. Similarly, if the existing technology requires 20\% (\resp 40\%) less power, then a membrane  perm-selectivity of at least 81 (\resp 231) is needed in order for the membrane cascade to operate with the same power as the existing technology. We note that energy is only one component of the overall separation cost. While we only discuss energy here, in practice, one will have to additionally account for the capital costs for a proper comparison.

%The above approach provides a better guidance than the upper bound curve on the Robeson plot, because it accounts for the energy requirement of alternate technology.
% allows comparison with alternate technology and also allows optimal choice of selectivity and permeability among membranes.
% helps experimental researchers in recognizing the target selectivity at which novel membrane materials will compete favorably with the existing alternative. As such, this targeting approach

\begin{table}[!ht]
     \centering
        \begin{tabular}{C{3cm} C{3cm}}
         \toprule
         p-xylene perm-selectivity &  Total power (kW)\\
         \midrule
         50   & 1,416 \\
         63   & 1,274 \\
         81   & 1,133 \\
         231 & 850 \\
         \bottomrule
        \end{tabular}
        \caption{Optimum power consumption for the cascade in Figure \ref{Fig_cascades_px_ox_4sgt} for different values of p-xylene perm-selectivity.}
        \label{Table_selectivities_px_ox}  
\end{table}

\subsubsection{Optimal four-stage cascade with two intermediate pumps}
Unlike compressors, pumps are relatively cheap and easy to operate. Thus, the number of pumps used for liquid separations may not be a major concern. Nevertheless, operating a cascade beyond a certain number of pumps may not be economically attractive. In such cases, one can append \eqref{const_comp1}-\eqref{const_comp4} to the MINLP and solve it to obtain attractive cascades with limited number of pumps.\\ 

\begin{figure}[h!]
    \centering
    \includegraphics[width=0.95\textwidth]{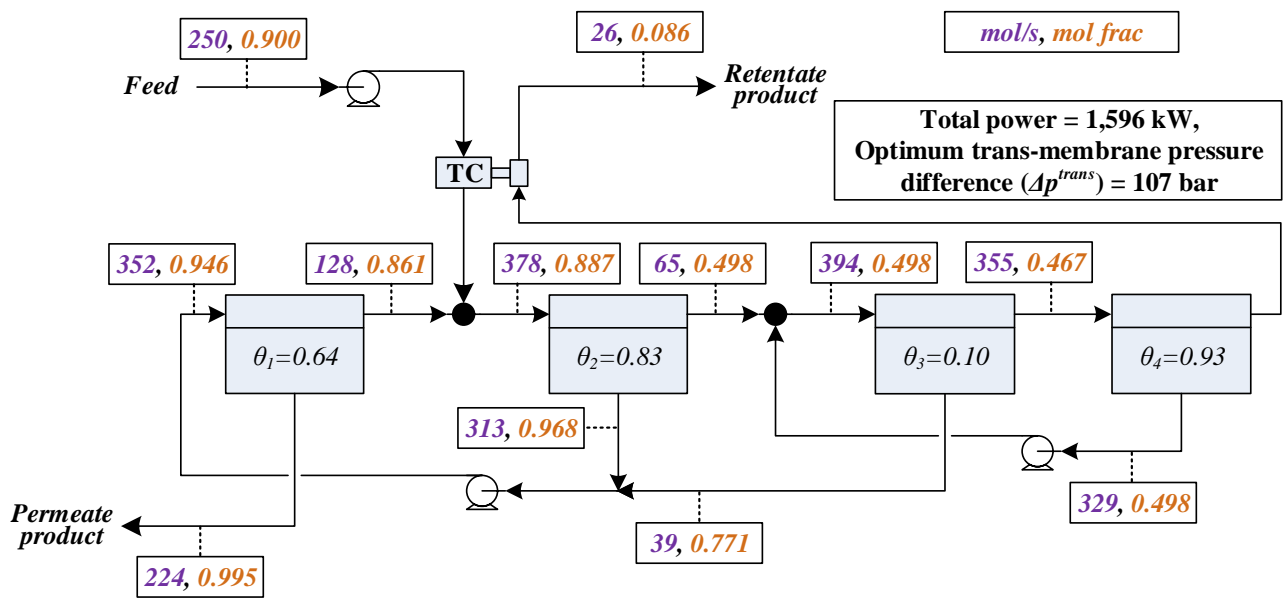}
    \caption{Optimal four-stage cascade with only two intermediate pumps for the separation of xylene isomers. All mole fractions shown in the box-legends correspond to p-xylene. Any discrepancy in the mass balances is due to the rounding of the flow and the composition values.}
    \label{Fig_cascades_px_ox_2}
\end{figure}

For the xylene mixture, we determine the optimal cascade requiring at most two intermediate pumps \ie{} $M_{comp}=2$. Solving (W) after appending \eqref{const_comp1}-\eqref{const_comp4} yields the solution shown in Figure \ref{Fig_cascades_px_ox_2}. Similar to Figure \ref{Fig_cascades_px_ox_4sgt}, the feed mixture is fed at the second stage in Figure \ref{Fig_cascades_px_ox_2}. However, in the latter arrangement, the permeate streams from both the second and third stages are mixed and recycled to the first stage. At optimum operation, this cascade requires 13\% more energy than the one in Figure \ref{Fig_cascades_px_ox_4sgt}. An economic analysis of these arrangements reveals whether the reduction in the capital cost due to elimination of a pump overweights the increase in operating cost due to the increase in power consumption or not. 
% is beneficial enough to build a cascade that consumes more energy. 

\subsubsection{Optimal three-stage cascade}
\begin{figure}[h!]
    \centering
    \includegraphics[width=1\textwidth]{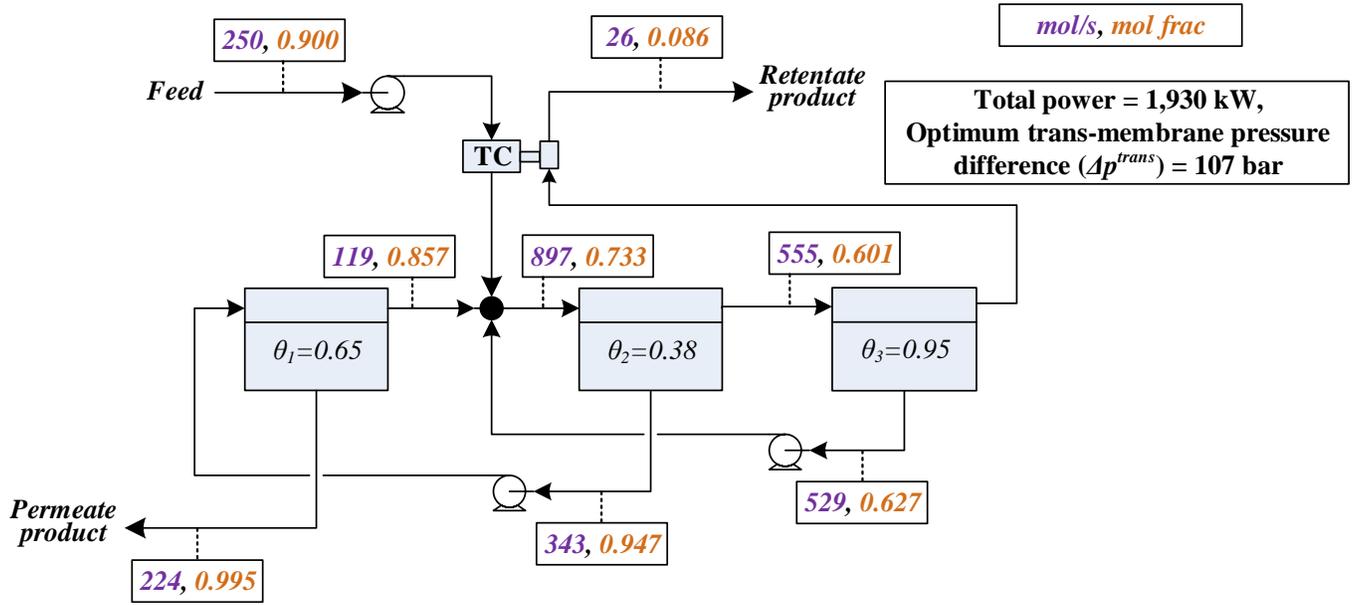}
    \caption{Optimal three-stage cascade for the separation of p-xylene. All mole fractions shown in the box-legends refer to p-xylene. Any discrepancy in the mass balances is due to the rounding of the flow and the composition values.}
    \label{Fig_cascades_px_ox_3stg}
\end{figure}

The main purpose of this subsection is to demonstrate that solving the MINLP with \eqref{const_comp1}-\eqref{const_comp4} often leads to better solutions than simply solving (W) with a fewer number of stages. For the xylene separation, a membrane cascade with only two intermediate pumps can also be designed by solving (W) with $N=3$. Figure \ref{Fig_cascades_px_ox_3stg} shows the optimal cascade along with its optimal operating condition obtained by solving (W) with $N=3$. Interestingly, the optimal three-stage cascade consumes significantly more energy than the cascade in Figure \ref{Fig_cascades_px_ox_2} even though both the arrangements require the same number of pumps. This example illustrates that \eqref{const_comp1}-\eqref{const_comp4} allow discovery of energy efficient cascades with limited number of pumps/compressors which are otherwise not easy to obtain directly by solving the MINLP without these constraints.

\section{Concluding remarks}

In this article, we addressed the problem of identifying energy efficient membrane cascades for the separation of a given binary liquid or gaseous mixture through dense membranes. The main highlights of the article are the following. First, to model the permeation process, we used a simple yet reasonably accurate model that does not require the solution of a detailed set of transport equations on either side of the membrane. This model was derived by solving a DAE system analytically. An advantage of the analytical solution is that it circumvents the need for discretization, where the latter technique requires a trade-off between model accuracy and complexity. Second, we proposed a novel MINLP that is applicable for both gaseous and liquid mixtures. Third, we demonstrated the need for a global optimality certificate on a numerical example, where we showed that local solution may consume significantly more energy than the globally optimal cascade. Unfortunately, state-of-the-art global solvers fail to solve the problem to a specified optimality tolerance in a reasonable amount of time. To aid the convergence of the solvers, we derived additional cuts that utilize physical intuition and the mathematical structure of the equations. Our computational experiments demonstrated that these cuts significantly improve the convergence of the global solvers. Finally, we used our model to find attractive cascades for two important industrial separations.

 \section*{Acknowledgements}
Jose Adrian Chavez Velasco thanks COMEXUS-Fulbright, and CONACYT for providing him with financial support. We thank Professor Ryan Lively for sharing data in reference \cite{koh2016reverse} with us.

\appendix

\section{Appendix}
\subsection{Reformulation of the flux equation}
\label{sec:flux-reformulation}
We begin with the definition of the local mole fraction on the permeate side
\begin{align*}
& y = \frac{n_A}{n_A + n_B} = \frac{S \left[ x - y \exp(-C_A u)\right]}{ S \left[ x - y \exp(-C_A u)\right] + \left[ (1-x) - (1-y) \exp(-C_B u)\right]}\\
\implies & y = \frac{S x - S y \exp(-C_A u)}{x(S-1) - y \left[S \exp(-C_A u) - \exp(-C_B u)\right] + 1-\exp(-C_B u)}
\end{align*}
Multiplying both sides with the denominator of the right hand side leads to 
\begin{align*}
& yx(S-1) + y = Sx - \left[S \exp(-C_A u) - \exp(-C_B u)\right]y (1-y)\\
\implies & \frac{y}{S-(S-1) y} = x - \frac{\left[S \exp(-C_A u) - \exp(-C_B u)\right]y (1-y)}{S-(S-1) y}\\
\implies & y-x = y - \frac{y}{S-(S-1) y}
- \frac{\left[S \exp(-C_A u) - \exp(-C_B u)\right]y (1-y)}{S-(S-1) y}\\
 \implies & y-x = \left[S (1-e^{-C_A u}) - (1-e^{-C_B u})\right] \frac{y(1-y)}{S-(S-1) y}
\end{align*}

we use partial fraction decomposition to rearrange the right hand side as follows
\begin{align*}
& y-x =  \left[\frac{S (1-e^{-C_A u}) - (1-e^{-C_B u})}{ (S-1)^2 }\right] \left(1 + (S-1) y - \frac{S}{S-(S-1) y} \right).
\end{align*}
Finally, we obtain \eqref{genfluxrel_2} by defining the term inside square brackets as $k$. 

\subsection{Derivation of the Analytical Solution}
\label{sec:analytic-der}
Consider the DAE system 
\begin{align*}
    & \frac{dx}{df} = \frac{y-x}{f}, \quad x(f^\feed) = x^\feed, \quad f \in [f^\feed,f^\ret]\\
    & y-x = k \left(1 + (S-1) y - \frac{S}{S-(S-1) y} \right)
\end{align*}
This implies
\begin{align*}
    & \int_{x^\feed}^{x^\ret} \frac{dx}{y-x} = \int_{x^\feed}^{x^\ret} \frac{d(y-(y-x))}{y-x}  = \int_{f^\feed}^{f^\ret} \frac{df}{f},\\
   \implies & \int_{y^\feed}^{y^\exit} \frac{dy}{y-x} - \int_{y^\feed-x^\feed}^{y^\exit-x^\ret} \frac{d(y-x)}{y-x} = \int_{f^\feed}^{f^\ret} \frac{df}{f}\\
   \implies & \frac{1}{k(S-1)^2}\int_{y^\feed}^{y^\exit} \frac{S-(S-1) y}{y(1-y)} dy - \int_{y^\feed-x^\feed}^{y^\exit-x^\ret} \frac{d(y-x)}{y-x} = \int_{f^\feed}^{f^\ret} \frac{df}{f} \\
   \implies & \frac{1}{k(S-1)^2}\int_{y^\feed}^{y^\exit} \left[\frac{S}{y}+\frac{1}{1-y}\right]dy - \ln \frac{y^\exit-x^\ret}{y^\feed-x^\feed} = \ln \frac{f^\ret}{f^\feed}\\
   \implies & \frac{1}{k(S-1)^2} \left[S \ln \frac{y^\exit}{y^\feed} - \ln \frac{1-y^\exit}{1-y^\feed}\right] - \ln \frac{y^\exit-x^\ret}{y^\feed-x^\feed} = \ln \frac{f^\ret}{f^\feed}
\end{align*}

% \pagebreak

\subsection{Derivation of $dy/dx$ used in Property P3}\label{derivation_dydx}
We differentiate \eqref{genfluxrel_2} with respect to $y$ to obtain

\begin{align*}
    \frac{dx}{dy} & =1-k (S-1) + \frac{k(S-1)S}{\left[ S- (S-1) y \right]^2}\\
        & = \frac{1}{y} \left\{ y- k(S-1)y + \frac{k y (S-1)S}{\left[ S- (S-1) y  \right]^2}   \right\}\\
        & \overset{\eqref{genfluxrel_2}}{=} \frac{1}{y} \left\{ x +k - \frac{k S }{S- (S-1) y}  + \frac{k y (S-1)S}{\left[ S- (S-1) y \right]^2}  \right\}  \\ 
        & = \frac{1}{y} \left\{ x + \frac{k (S-1)^2 y^2}{\left[ S (1-y) + y \right]^2}  \right\}. \label{dydx}
    \end{align*}
Therefore, 
\begin{align*}
    \frac{dy}{dx} = \left( \frac{dx}{dy} \right)^{-1} = y \left[ x + \frac{k (S-1)^2 y^2}{\left[ S (1-y) + y \right]^2} \right]^{-1}.
\end{align*}

\begin{comment}
\begin{subequations}
    \begin{align}
        & \frac{dx}{dy} = \frac{1}{y} \left\{ y+ y k(1-S) + \frac{k y (S-1)S}{\left[ S (1-y) + y  \right]^2}   \right\}\\
        & \frac{dx}{dy} = \frac{1}{y} \left\{ x- \left[y-\frac{k (S-1)^2 (1-y) y}{S (1-y) +y} \right] + y + y k(1-S) + \frac{k y (S-1)S}{\left[ S (1-y) + y \right]^2}  \right\}  \\ 
        & \frac{dy}{dx} = \frac{1}{y} \left\{x + \frac{k (S-1)^2 (1-y) y}{S (1-y) + y} + y k (1-S) + \frac{k y (S-1) S}{\left[ S (1-y) + y \right]^2}  \right\}\\
        & \frac{dx}{dy}= \frac{1}{y} \left \{x + \frac{k y^2 -2 k S y^2 + k S^2 y^2}{\left[ S (1-y) + y \right]^2}  \right\}\\ 
        & \frac{dy}{dx} = \frac{1}{y} \left\{ x + \frac{k (S-1)^2 y^2}{\left[ S (1-y) + y \right]^2}  \right\} \label{dydx}
    \end{align}
\end{subequations}

Since all the terms on the right hand side of \eqref{dydx} are non-negative, $\frac{dx}{dy}\geq 0$, and $\frac{dy}{dx} = \left( \frac{dx}{dy} \right)^{-1}  \ge 0$. This implies that $y$ increases monotonically with $x$.
\end{comment}

\bibliography{bibliography} 
\bibliographystyle{ieeetr}

\end{document}